\documentclass{article}
\usepackage{iclr2027_conference,times}
\usepackage[T1]{fontenc}
\usepackage{amsmath,amssymb}
\usepackage{graphicx}
\usepackage{algorithm}
\usepackage{algpseudocode}
\usepackage[hidelinks]{hyperref}
\usepackage{url}
\usepackage{booktabs}
\usepackage{placeins}
\usepackage{enumitem}
\usepackage{amsthm}
\usepackage{xcolor}
\usepackage[most]{tcolorbox}
\usepackage{adjustbox}
\usepackage{subcaption}
\usepackage{wrapfig}
\usepackage{titletoc}
\usepackage{needspace}
\usepackage{ragged2e}
\usepackage{array,tabularx}
\usepackage{capt-of}

\usepackage{amsmath,amsfonts,bm}

\def\eqref#1{equation~\ref{#1}}

\def\1{\bm{1}}

\DeclareMathAlphabet{\mathsfit}{\encodingdefault}{\sfdefault}{m}{sl}
\SetMathAlphabet{\mathsfit}{bold}{\encodingdefault}{\sfdefault}{bx}{n}

\newcommand{\ph}[1]{%
  \textcolor{red!75!black}{\texttt{\{#1\}}}%
}

\newtcolorbox{promptbox}[2]{
  enhanced,
  breakable,
  colback=white,
  colframe=black,
  colbacktitle=black,
  coltitle=white,
  title={#1},
  fonttitle=\bfseries,
  boxrule=0.9pt,
  arc=1.5mm,
  outer arc=1.5mm,
  left=8pt,
  right=8pt,
  top=7pt,
  bottom=7pt,
  before skip=8pt,
  after skip=10pt,
  label={box:#2}
}

\newlist{promptlist}{description}{3}
\setlist[promptlist]{
  leftmargin=1.8em,
  labelsep=0.45em,
  itemsep=3pt,
  topsep=3pt,
  parsep=0pt,
  partopsep=0pt
}
\newcommand{\methodname}{DiSH}

\title{
Harnessing LLMs with Direction-Scope Control for Combinatorial Solver Refinement
}

\author{
Yili Cao$^{1}$
\quad
Anqi Li$^{1,\dagger}$
\quad
Haoran Li$^{3}$
\quad
Yongtang Shi$^{1,2}$
\quad
Zicheng Zhang$^{4,\dagger}$ 
\thanks{The authors are listed in alphabetical order.}\\
$^{1}$ School of Mathematical Sciences, Nankai University 
\\
$^{2}$ Center for Combinatorics and LPMC, Nankai University
\\
$^{3}$ David R. Cheriton School of Computer Science, University of Waterloo
\\
$^{4}$ JD.com \quad
$^{\dagger}$ Corresponding author
}

\iclrfinalcopy

\begin{document}

\maketitle

\begin{abstract}
Harnesses are increasingly unlocking the potential of large language model (LLM) agents by making their actions more controllable. In automatic heuristic design (AHD) for combinatorial optimization (CO), LLMs iteratively revise solvers using execution and verification feedback. However, these revisions often break previously validated solver behavior, making refinement unstable and costly. We characterize this pattern as \emph{feasibility oscillation} and document it across CO benchmarks, AHD methods, and backbone LLMs. Tracing this oscillation, our analysis reveals that feasibility degradation often follows large solver revisions, while substantial feedback redundancy makes it harder to identify which failures deserve attention. Motivated by these findings, we introduce DiSH, a harness that jointly controls revision direction and scope. It selects representative failure evidence to guide each update and progressively narrows the revision scope as feasible behavior accumulates. Across nine CO benchmarks, DiSH consistently improves QYI scores and reduces feasibility oscillation, while roughly halving feedback token consumption relative to prior methods that use full feedback. These results establish DiSH as a strong baseline for solver refinement under harness control.

\end{abstract}

\section{Introduction}
\begin{wrapfigure}{r}{0.47 \textwidth}
    \centering
    \vspace{-1.6em}
    \includegraphics[width=\linewidth]{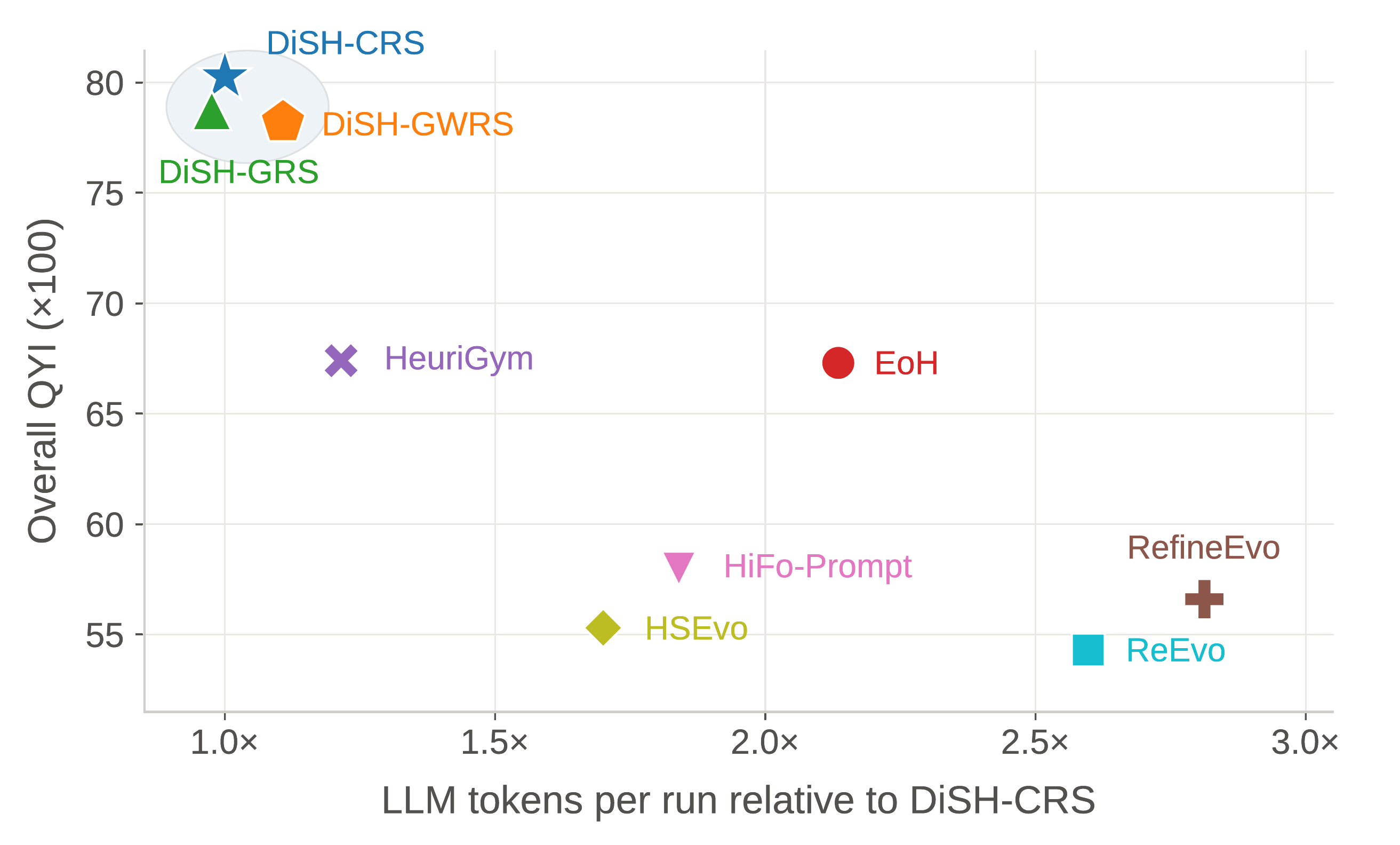}
    \vspace{-2em}
    \caption{Overall QYI versus token cost. }
    \label{fig:cost_tradeoff1}
    \vspace{-1em}
\end{wrapfigure}

Automatic heuristic design seeks to automate development for combinatorial optimization algorithms. Traditional methods search over heuristics predefined by human experts \citep{burke2013,pillay2018}, while large language models enable solvers to be generated and revised through iterative feedback \citep{funsearch,eoh,heurigym}. This execution-and-feedback loop naturally aligns with the harness, which has become a foundation for building LLM agents \citep{coala,sweagent,openhands}. Generally, a harness is the system surrounding an LLM that structures its observations and actions, executes its decisions, and returns results. In LLM-based AHD, it therefore determines how execution and verification evidence informs subsequent solver revisions.

Despite this central role, the harness-like components of existing LLM-based AHD methods remain largely limited to executing and verifying candidate solvers and returning evaluation feedback to the LLM \citep{funsearch,eoh,reevo,heurigym}. This  is particularly problematic for CO solvers, whose constraints, decisions, and long-horizon logic are tightly coupled. A reported execution error or constraint violation identifies the symptom of a failure, but rarely specifies its root cause in the program or the extent of revision required to correct it. The LLM therefore decides which failures deserve attention and how extensively the solver should be modified. Existing methods expose what went wrong, but leave how to revise largely implicit in model generation. Recent methods regulate the search more finely, but they typically couple revision scope to search operators rather than constraining it independently before generation~\citep{mctsahd,a2dept,hifoprompt,refineevo}, while revision priorities largely depend on the model's interpretation of richer feedback~\citep{lago,revel,pathwise,ahdagent}. We discuss these lines of work in detail in Appendix~\ref{sec:related}.

This leaves LLM-based AHD with a central challenge: \textit{iterative revision does not reliably and efficiently accumulate improvement across rounds.} A revision may repair current failures while disrupting verified behavior, forcing later rounds to recover lost progress. We characterize this repeated loss and recovery as \emph{feasibility oscillation}, and find it consistently across CO benchmarks, AHD methods, and backbone LLMs, indicating that it is a systemic rather than an isolated issue.  Our analysis links feasibility loss to large changes between consecutive solvers and reveals substantial redundancy in the feedback guiding those changes, which obscures repair priorities. These findings motivate joint control of revision scope, which limits how extensively the solver may change, and revision direction, which determines which failures guide the update.

\begin{figure}[t]
    \centering
    \vspace{-2em}
    \includegraphics[width=1.0\linewidth]{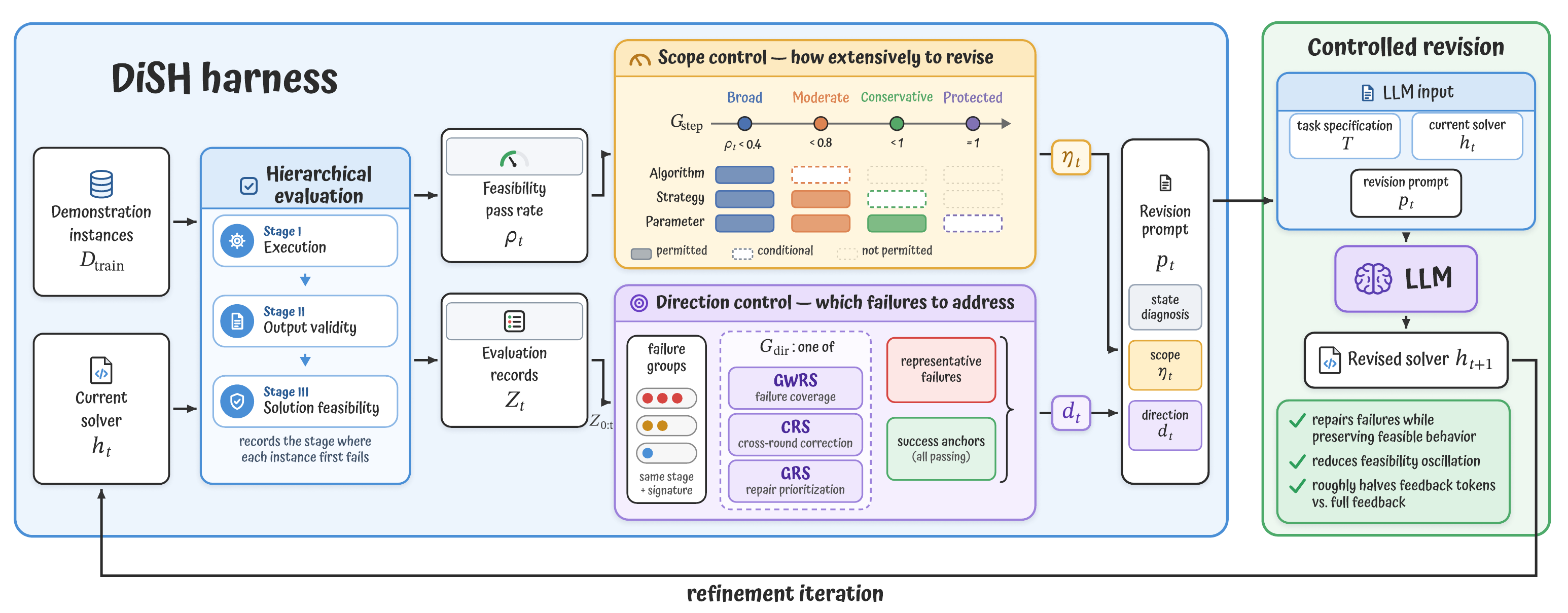}
    \caption{ Overview of the DiSH framework. DiSH converts evaluation feedback into two complementary controls over what to revise and how extensively to revise it. Scope control adjusts the permitted revision magnitude, while direction control selects representative failures and successful instances to focus the revision. Together, they form the revision prompt that guides the LLM to repair remaining failures while preserving previously feasible behavior across iterations. }
    \label{fig:overview}
\end{figure}

Guided by this analysis, we introduce DiSH, a harness framework designed to mitigate feasibility oscillation through joint control of revision direction and scope.  As illustrated in Figure~\ref{fig:overview}, DiSH first organizes execution and verification outcomes through hierarchical evaluation, preserving both instance-level failure evidence and the aggregate feasibility state. Direction control distills redundant records into representative repair evidence, making the target of each revision explicit. Scope control progressively restricts the permitted modification as feasible behavior accumulates, reducing the risk that a new revision disrupts verified behavior. Together, these components turn raw evaluation outcomes into revision instructions, guiding iterative revision toward cumulative refinement.

Figure~\ref{fig:cost_tradeoff1} summarizes the performance and efficiency comparison. Across nine CO benchmarks, the three DiSH variants achieve higher average QYI while consuming fewer LLM tokens than existing AHD methods. Token-accounting details and the per-backbone search-cost breakdown are provided in Appendix~\ref{app:search_cost}. DiSH also consistently reduces feasibility oscillation and roughly halves feedback-token usage relative to full feedback. Further analysis shows that narrower scope instructions lead to smaller realized revisions, while compact direction evidence preserves repair effectiveness. These results establish DiSH as a strong baseline for harness-controlled solver refinement.

Our main contributions are summarized as follows:
\begin{itemize}[
    topsep=0pt,
    partopsep=0pt,
    itemsep=2pt,
    parsep=0pt,
    leftmargin=0pt
]
    \item We identify {feasibility oscillation} as a central obstacle to cumulative solver refinement and systematically characterize it across CO benchmarks, AHD methods, and backbone LLMs. Our analysis associates feasibility loss with large solver revisions and reveals substantial redundancy in the feedback used to guide these revisions.

    \item \textit{We introduce DiSH, the first AHD harness to our knowledge to demonstrate stable and cumulative improvement across revisions}.  DiSH jointly controls revision direction and scope by distilling repeated failures into representative repair evidence and progressively restricting solver modifications as feasible behavior accumulates.

    \item We evaluate DiSH across nine CO benchmarks and show consistent improvements in QYI and refinement stability. DiSH reduces feasibility oscillation and roughly halves feedback-token usage relative to full-feedback refinement, while preserving repair effectiveness and generalizing across direction strategies and backbone LLMs.
\end{itemize}

\enlargethispage{2\baselineskip}

\section{Problem Formulation}\label{subsec:overview}

LLM-based automatic heuristic design (AHD) commonly improves solvers through iterative code revision~\citep{eoh,reevo,heurigym,refineevo}. Given a task specification $\mathcal{T}$ and demonstration instances $\mathcal{D}_{\mathrm{train}}$, the current solver $h_t$ is executed and evaluated on these instances. The resulting feedback $\mathrm{F}_t$ guides the next revision:
\begin{equation}
h_{t+1} \sim \mathcal{M}_{\mathrm{LLM}} \left( \cdot \mid \mathcal{T},h_t,\mathrm{F}_t \right),
\label{eq:feedback_refinement}
\end{equation}
where $\mathcal{M}_{\mathrm{LLM}}$ denotes the conditional distribution over solver programs induced by the LLM. Applying~(\ref{eq:feedback_refinement}) for $K$ steps yields a trajectory $h_0,h_1,\ldots,h_K$. Define the feasibility pass rate of $h_t$ as
\begin{equation}
\rho_t = \frac{1}{|\mathcal{D}_{\mathrm{train}}|} \sum_{x\in\mathcal{D}_{\mathrm{train}}} \mathbb{I}\!\left[ h_t \text{ produces a feasible solution for } x \right].
\label{eq:feasibility_rate}
\end{equation}
\begin{wrapfigure}{r}{0.46\linewidth}
    \centering
    \vspace{-2em}
    \includegraphics[width=\linewidth]{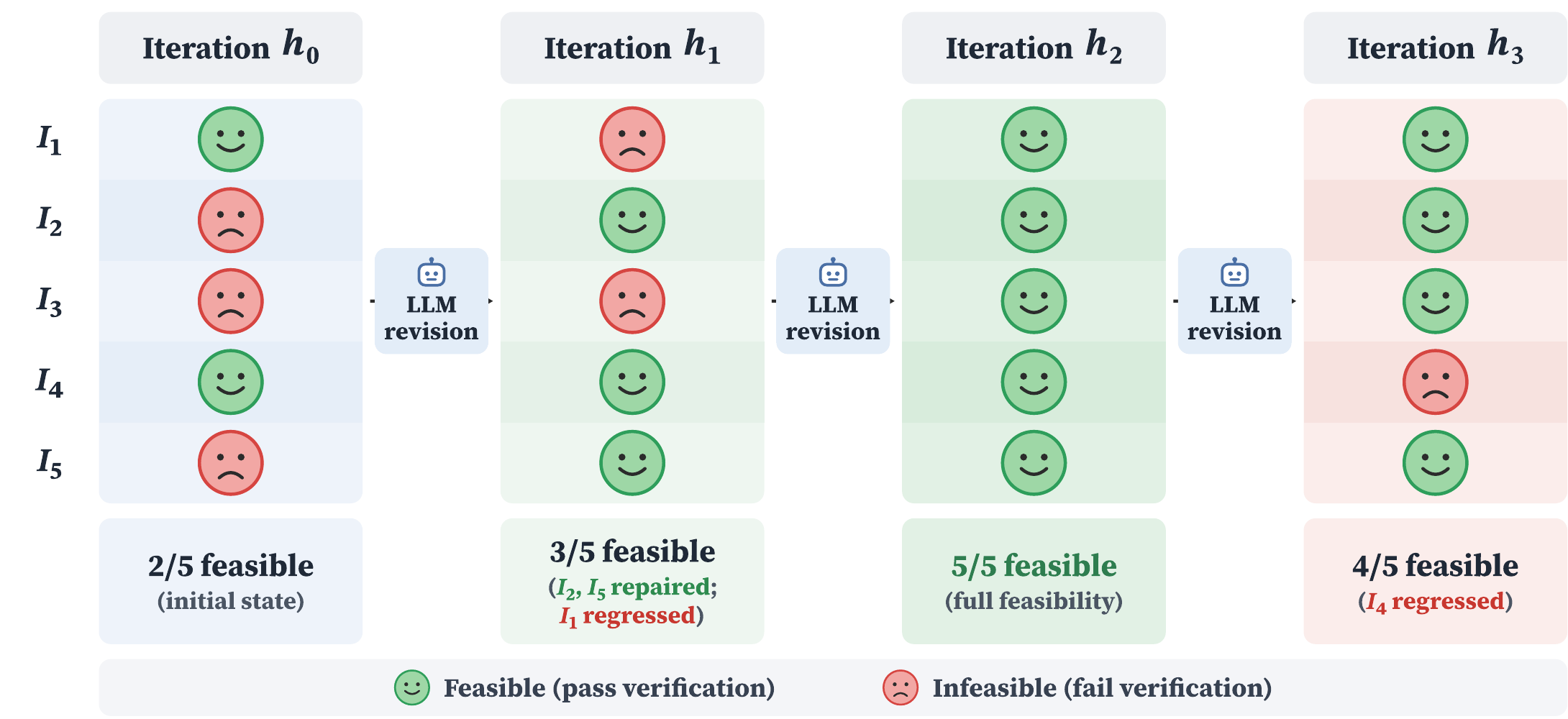}
    \vspace{-2em}
    \caption{Instance-level feasibility oscillation.}
    \label{fig:instance_feasibility_oscillation}
    \vspace{-1em}
\end{wrapfigure}
\textbf{Feasibility Oscillation.} A revision may repair current failures while making previously feasible instances fail. Repeated losses and recoveries of feasibility lead to \emph{feasibility oscillation}. A decrease in $\rho_t$ indicates a net loss of feasibility, but individual regressions can also occur when the pass rate increases. Figure~\ref{fig:instance_feasibility_oscillation} illustrates this distinction. The first revision repairs $I_2$ and $I_5$ but makes $I_1$ infeasible, raising the pass rate from $2/5$ to $3/5$ despite this regression. The next revision restores $I_1$ and repairs $I_3$, reaching full feasibility. A further revision, however, makes $I_4$ infeasible and lowers the pass rate to $4/5$. Such regressions consume additional revision budget and hinder stable progress, motivating a closer analysis of the underlying solver revisions.

\textbf{Revision Distance as a Lens on Feasibility Oscillation.} We measure the distance between consecutive solvers to examine how the revisions relate to feasibility oscillation. Let $T(h)$ denote the set of contiguous token sequences extracted from solver $h$. We use their Jaccard distance as an empirical measure of revision distance (Appendix~\ref{app:revision_distance}):
\begin{equation}
\Delta(h_{t+1},h_t) = 1- |T(h_t)\cap T(h_{t+1})|/|T(h_t)\cup T(h_{t+1})|.
\label{eq:revision_magnitude}
\end{equation}
\begin{wrapfigure}{r}{0.46\linewidth}
    \centering
    \vspace{-1em}
    \includegraphics[width=\linewidth]{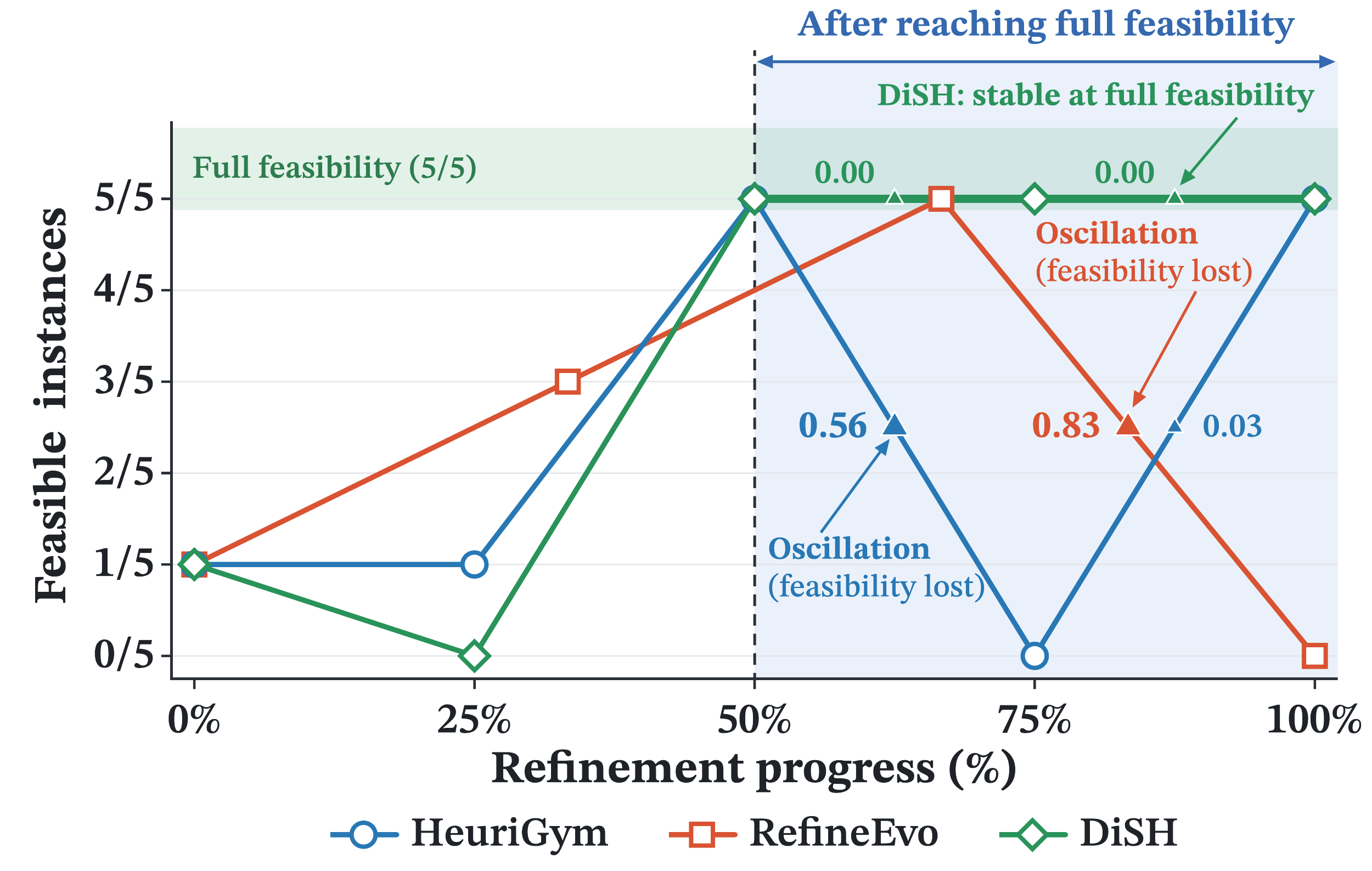}
    \vspace{-2em}
    \caption{Revision distance and oscillation.}
    \label{fig:probing_oscillation}
    \vspace{-1em}
\end{wrapfigure}
Smaller distances indicate greater code retention, while larger distances indicate more extensive changes. Figure~\ref{fig:probing_oscillation} compares trajectories from HeuriGym~\citep{heurigym}, RefineEvo~\citep{refineevo}, and \methodname{}. The curves show feasible instance counts, and triangles show revision distances $\Delta(h_{t+1},h_t)$. HeuriGym and RefineEvo drop from full feasibility to $0/5$ after substantial revisions with distances of $0.56$ and $0.83$, respectively. In contrast, \methodname{} maintains full feasibility once reached, with subsequent distances reported as $0.00$. These examples suggest that extensive revisions can significantly disrupt established feasibility, leading to unexpected oscillation.

\section{\methodname{}: Direction-and-Scope Harness for CO Solver Refinement}
\label{sec:method}

Combinatorial optimization (CO) solvers often involve sequences of interdependent decisions. As discussed in Section~\ref{subsec:overview}, revising their code can disrupt feasible behavior instead of improving it. We introduce \methodname{}, a harness that explicitly controls what to revise and how extensively to change the solver. It aims to stabilize feasibility gains while using feedback tokens efficiently. We first describe the controlled revision loop (Section~\ref{subsec:harness_control}). We then present hierarchical evaluation (Section~\ref{subsec:feasibility_evaluation}), adaptive revision scope (Section~\ref{subsec:scope}), and revision-direction construction (Section~\ref{subsec:direction}).

\subsection{From Iterative Revision to Harness-Controlled Solver Refinement}
\label{subsec:harness_control}

We turn evaluation feedback into explicit controls over the direction and scope of solver revisions.

\textbf{Why Revisions Need Control.} Feasibility oscillation highlights that two issues limit stable refinement. First, a single LLM call may not resolve failures involving closely linked constraints, decisions, and solver components. Second, repeating revisions can also undo earlier repairs. Thus, \emph{revision} changes the solver, whereas \emph{refinement} should accumulate feasible behavior across rounds.

Given $h_t$, \methodname{} produces per-instance evaluation records $\mathcal{Z}_t$ and the feasibility pass rate $\rho_t$. It then constructs a revision direction $d_t$, a scope instruction $\eta_t$, and a prompt $p_t$ for the next LLM revision:
\begin{equation}
\begin{aligned}
\text{Evaluation: } & (\mathcal{Z}_t,\rho_t) = \mathcal{E}(h_t,\mathcal{D}_{\mathrm{train}}), \\
\text{Construction: } & d_t = \mathcal{G}_{\mathrm{dir}}(\mathcal{Z}_{0:t}), \quad \eta_t = \mathcal{G}_{\mathrm{step}}(\rho_t), \\
\text{Generation: } & p_t = \mathcal{H}(\mathcal{Z}_t,d_t,\eta_t), \quad h_{t+1} \sim \mathcal{M}_{\mathrm{LLM}} (\cdot\mid\mathcal{T},h_t,p_t).
\end{aligned}
\label{eq:controlled_refinement}
\end{equation}
Here, $\mathcal{E}$ evaluates instances in successive stages. $\mathcal{G}_{\mathrm{dir}}$ selects and organizes revision targets from current and historical evaluation evidence, while $\mathcal{G}_{\mathrm{step}}$ determines the permitted scope from $\rho_t$. $\mathcal{H}$ combines the evaluation records and these controls into revision instructions $p_t$, which is then used to generate the next solver $h_{t+1}$. Appendix~\ref{app:prompts} provides the complete prompt structure and examples.

Rather than exposing all observed failure records directly to the LLM, direction control prioritizes them into representative repair evidence while retaining previously verified behavior. Scope control permits less extensive changes at higher pass rates, encouraging smaller revision distances $\Delta(h_{t+1},h_t)$. Together, these controls aim to reduce regressions and support stable feasibility gains within the finite revision budget. Algorithm~\ref{alg:dish_refinement} in Appendix~\ref{app:harness-loop} gives the complete refinement loop.

\subsection{Hierarchical Feasibility Evaluation and Feedback}
\label{subsec:feasibility_evaluation}

We evaluate each instance in three stages to identify where it fails and provide feedback for revision.

\textbf{From Outcome Evaluation to Failure Localization.} Existing AHD pipelines mainly use execution and verification outcomes to score solver candidates~\citep{funsearch,eoh,reevo} or provide feedback~\citep{heurigym}, without explicitly localizing failures for revision control. However, execution errors, invalid outputs, and constraint violations arise from various solving stages and require different repairs. \methodname{} therefore constructs a hierarchical evaluation pipeline that records both the failure stage and the associated diagnostic information.

\textbf{Stage-Wise Instance Evaluation.} To make evaluation feedback more diagnostic, \methodname{} applies three checks in order and records the first failed stage rather than reducing its outcome to a single score. For each instance $x\in\mathcal{D}_{\mathrm{train}}$, Stage~I checks execution within time and memory limits. Stage~II checks output validity, and Stage~III checks satisfaction of problem constraints. For example, a routing solver may time out on one instance but return a correctly formatted route that serves a customer twice on another. These failures call for different repairs: improving execution efficiency and correcting solution construction, respectively. Appendix~\ref{app:instance_evaluation} provides the detailed criteria.

\textbf{Hierarchical Evaluation Records.} The resulting record is $z_t(x)=(s_t(x),e_t(x),q_t(x))$, where $s_t(x)$ denotes the terminal evaluation stage, $e_t(x)$ contains the corresponding error information, and $q_t(x)$ is the objective value of a feasible solution. Depending on the terminal stage, $e_t(x)$ may contain an exception, a parsing error, or a constraint violation. Failure records have no objective value, while successful records have no error information. The records from all instances are collected as $\mathcal{Z}_t=\{z_t(x)\}_{x\in\mathcal{D}_{\mathrm{train}}}$. The hierarchical records support revision control at two levels. Aggregated across instances, they yield the feasibility pass rate $\rho_t$ defined in Eq.~(\ref{eq:feasibility_rate}), which informs the permitted revision scope. At the instance level, failure stages and diagnostics provide repair evidence, while successful records identify behavior to preserve.

\subsection{Feasibility-Adaptive Revision Scope}
\label{subsec:scope}

We adapt the permitted revision scope to the current pass rate, allowing less extensive changes as more instances become feasible.

\textbf{Conditioning Scope on Feasibility.} Failure evidence indicates what requires repair, but not how broadly the solver should be revised. At low $\rho_t$, fundamental failures may require substantial revisions. At high $\rho_t$, such revisions risk disrupting verified behavior, as shown in Figure~\ref{fig:probing_oscillation}. We therefore use $\eta_t=\mathcal{G}_{\mathrm{step}}(\rho_t)$ to progressively narrow the permitted revision scope as feasibility increases.

\textbf{Semantic Revision Levels and Modes.} To encourage smaller revision distances as feasibility improves, DiSH uses a prompt instruction $\eta_t$ to specify the permitted semantic scope. We distinguish three revision levels. \emph{Parameter} revisions adjust numerical constants without changing solver logic. \emph{Strategy} revisions modify decision rules or local procedures. \emph{Algorithm} revisions replace the core solving approach. The scope rule $\eta_t=\mathcal{G}_{\mathrm{step}}(\rho_t)$ selects among four modes with progressively tighter restrictions. \textsc{Broad} permits all three revision levels. \textsc{Moderate} favors strategy revisions and permits algorithm replacement only when failures indicate a structural mismatch. \textsc{Conservative} allows parameter changes and narrowly targeted strategy revisions but prohibits algorithm replacement. At full feasibility, \textsc{Protected} instructs the LLM to preserve the solver by default, permitting only beneficial parameter changes that do not compromise currently feasible instances. Thus, the permitted scope narrows as $\rho_t$ increases. Appendix~\ref{app:scope_rules} details the mode thresholds and permitted revision levels, and Appendix~\ref{app:stepsize} provides the corresponding prompts.

\needspace{1\baselineskip}
\begin{wrapfigure}{r}{0.4\linewidth}
    \centering
    \vspace{-1em}
    \includegraphics[width=\linewidth]{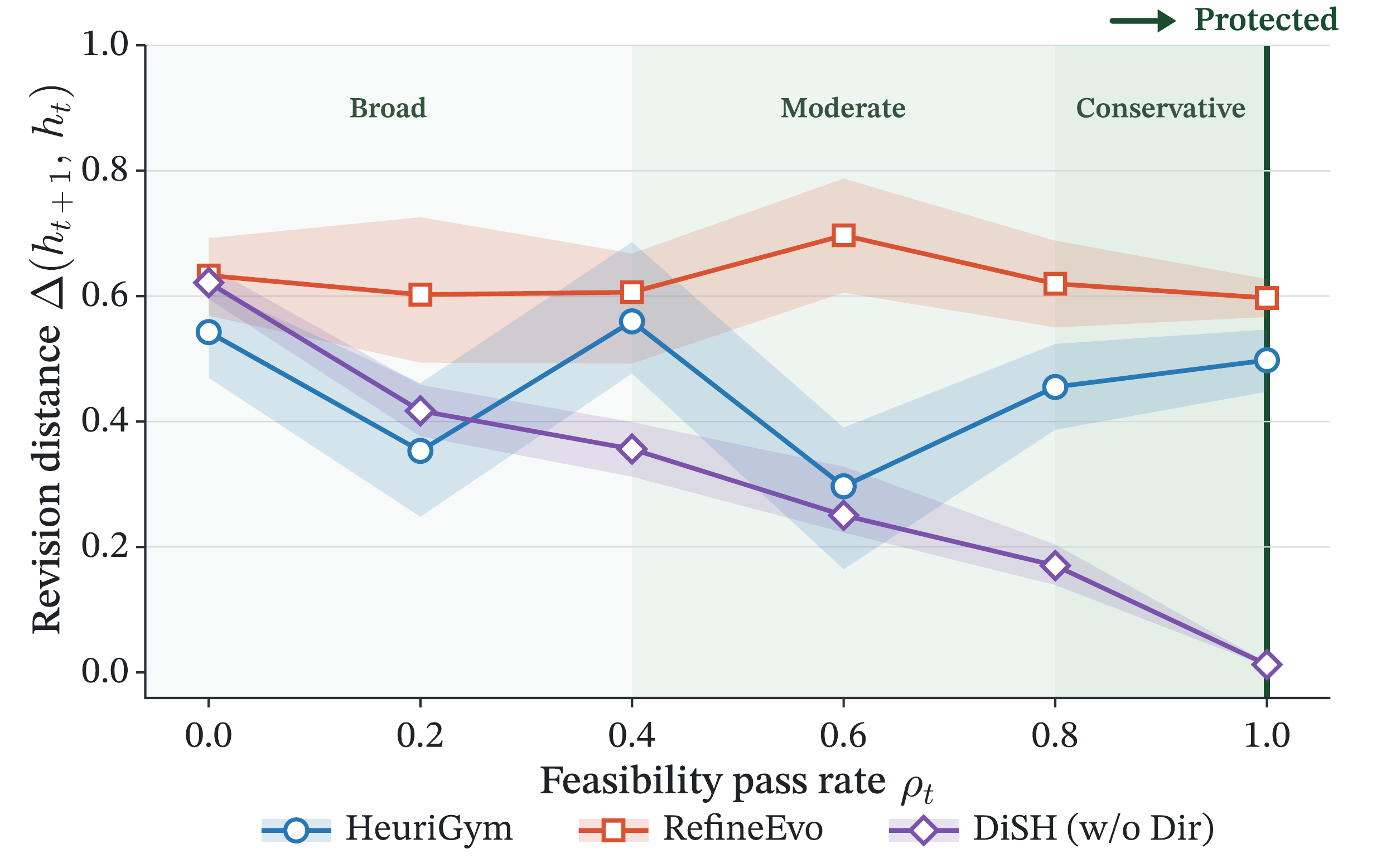}
    \vspace{-2em}
    \caption{Revision distance vs.\ $\rho_t$.}
    \label{fig:scope_preservation}
    \vspace{-1em}
\end{wrapfigure}
\textbf{Realized Revision Distance.} To examine whether scope instructions lead to corresponding changes in realized solver modification, we provide \emph{DiSH (w/o Dir)} as a controlled comparison, which retains full feedback and enables scope control. In Figure~\ref{fig:scope_preservation}, lines show mean revision distances, shaded ribbons indicate 95\% bootstrap confidence intervals, and background bands denote scope modes. As $\rho_t$ increases, the realized revision distance of \emph{DiSH (w/o Dir)} progressively contracts, whereas HeuriGym and RefineEvo do not exhibit the same trend. This pattern indicates that the revision scope specified by $\mathcal{G}_{\mathrm{step}}$ before generation is reflected in progressively smaller realized revisions as feasibility increases.

\subsection{Revision-Direction Construction}
\label{subsec:direction}

We build revision directions from representative failures and successful records, reducing repeated feedback while preserving useful repair guidance.

\begin{wrapfigure}{r}{0.4\linewidth}
    \centering
    \vspace{-1.5em}
    \includegraphics[width=\linewidth]{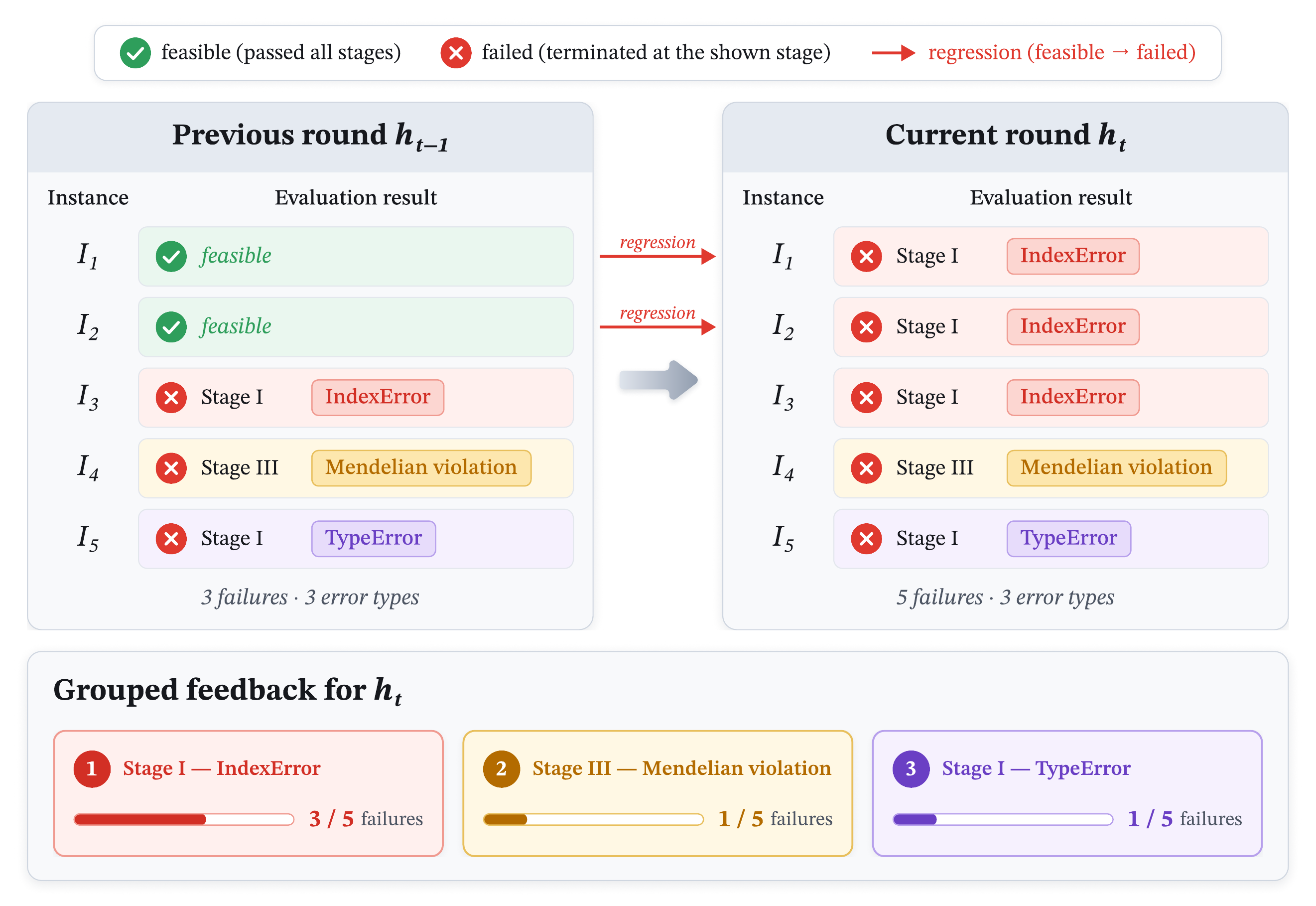}
    \vspace{-2em}
    \caption{Failure redundancy.}
    \label{fig:failure_redundancy_phenomenon}
    \vspace{-1.5em}
\end{wrapfigure}
\textbf{Distilling Evaluation Feedback into  Revision Direction.} Revision scope limits how extensively the solver may change, but a successful revision also depends on how the LLM uses evaluation feedback to decide what to change. The main obstacle is redundancy: many instance-level records describe the same failure, consuming tokens and making distinct repair needs harder to identify. Figure~\ref{fig:failure_redundancy_phenomenon} illustrates this redundancy in a refinement example. Previously feasible instances $I_1$ and $I_2$ now fail with the same \texttt{IndexError} as $I_3$. The number of failure records increases from three to five, while the number of distinct failure groups remains three. To further quantify feedback redundancy, we analyze 69,890 failure records collected across revision rounds from five methods on nine benchmarks. Of these records, 95.6\% share a failure pattern with another instance in the same evaluation. In 88.7\% of solver evaluations, at least two instances exhibit the same failure pattern. These results motivate grouping related failures and distilling representative evidence instead of passing all records to the LLM.

\textbf{Shared Failure-Group Representation.} Let $\mathcal{P}_t$ contain the instances that pass all evaluation stages and let $\mathcal{F}_t=\mathcal{D}_{\mathrm{train}}\setminus\mathcal{P}_t$ denote the set of failed instances. We group failed instances by their terminal stage and normalized error signature:
\begin{equation}
\mathcal{G}_t = \{G_t^{s,a}:G_t^{s,a}\neq\varnothing\}, \qquad \text{where } G_t^{s,a} = \{x\in\mathcal{F}_t:s_t(x)=s,\ \tau(e_t(x))=a\}.
\end{equation}
Here, $\tau$ maps the diagnostic information $e_t(x)$ to an error category $a$. Each group $G_t^{s,a}$ contains instances with the same failure stage $s$ and error category $a$, and $\mathcal{G}_t$ collects all nonempty failure groups. Across the collected 69,890 failure records, normalization yields 53 distinct error signatures $a$. Figures~\ref{fig:failure_distribution} and~\ref{fig:failure_cloud} in Appendix~\ref{app:failure_signatures} show their frequency distribution and representative normalized patterns, respectively. Instances within the same group therefore exhibit similar failure patterns and are expected to require similar revisions. A direction strategy selects a representative set $\mathcal{R}_t\subseteq\mathcal{F}_t$ from these groups. The records of $\mathcal{R}_t$, together with all successful records in $\mathcal{P}_t$, form the main content of $d_t$. We instantiate revision direction from three alternative perspectives that emphasize failure coverage, cross-round correction, and repair prioritization, respectively.

\begin{wrapfigure}{r}{0.31\linewidth}
    \centering
    \vspace{-1.8em}
    \includegraphics[width=\linewidth]{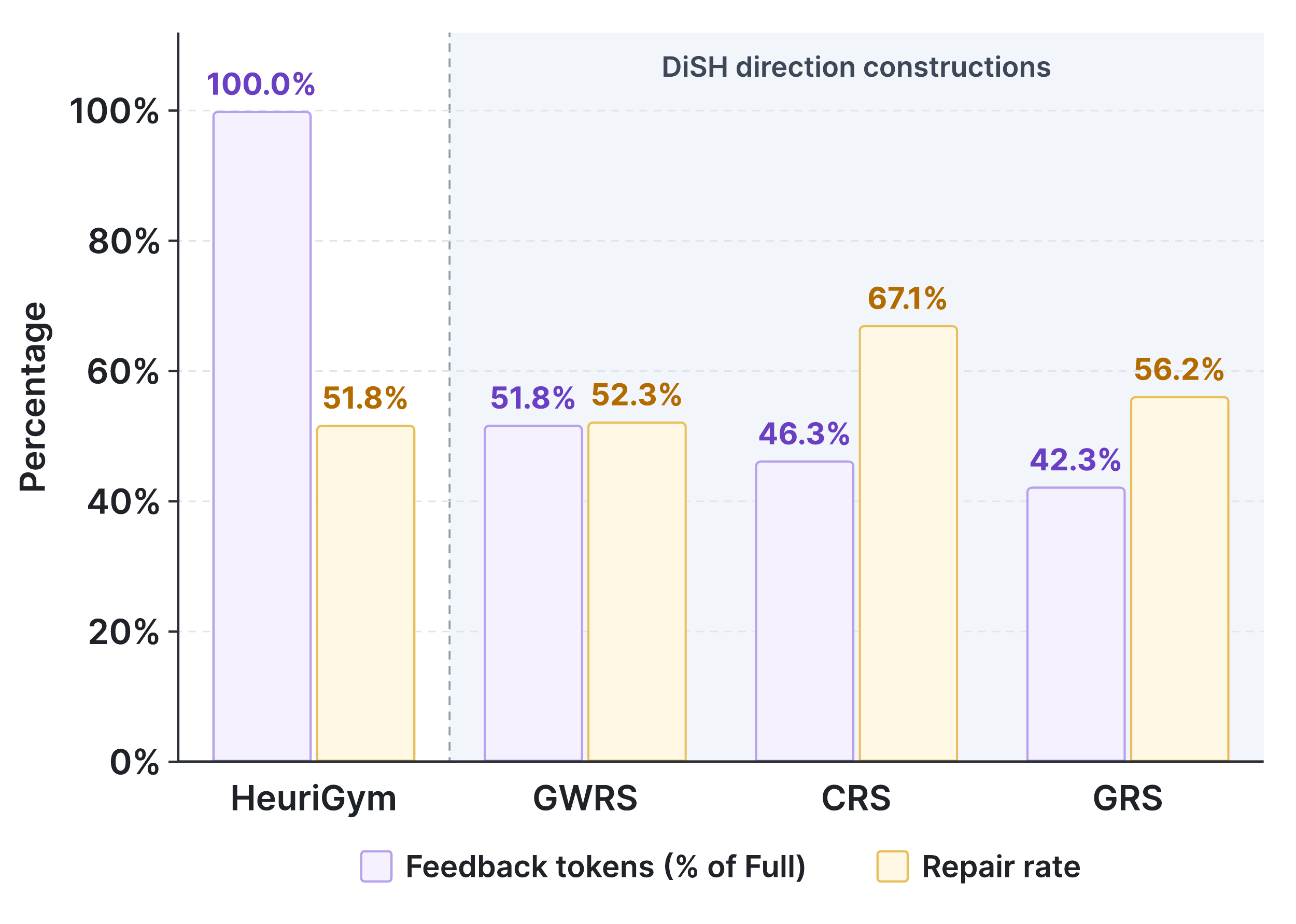}
    \vspace{-2em}
    \caption{Feedback efficiency.}
    \label{fig:feedback_efficiency_evaluation}
    \vspace{-1.5em}
\end{wrapfigure}
\textbf{Alternative Direction Constructions.} The direction constructor is chosen from three alternatives: $\mathcal{G}_{\mathrm{dir}}\in\{\mathcal{G}_{\mathrm{GWRS}},\mathcal{G}_{\mathrm{CRS}},\mathcal{G}_{\mathrm{GRS}}\}$. Each independently constructs $d_t$ from representative failures and successful records. As shown in Figure~\ref{fig:feedback_efficiency_evaluation}, they use 42.3--51.8\% of the tokens required by full feedback while preserving the measured rate at which failure modes are repaired in the next revision. Thus, selecting one direction strategy is sufficient to substantially reduce feedback without diminishing the measured repair effectiveness.
\begin{itemize}[leftmargin=0em, topsep=0.4em, itemsep=0.5em, parsep=0pt]
\item \textbf{Group-Wise Representative Sampling (GWRS)} selects one representative from every current failure group. For each $G\in\mathcal{G}_t$, it chooses the largest instance, $r_G=\operatorname*{arg\,max}_{x\in G}n(x)$, where $n(x)$ is the input size, with ties broken by instance identifier. This rule preserves group coverage without repeating similar errors. The largest instance is chosen to expose the failure on a structurally richer input. Groups are reconstructed after each evaluation: repaired instances provide examples of behavior to preserve, while unresolved failures remain available for later revisions.

\item \textbf{Correction Representative Sampling (CRS)} uses consecutive evaluations to relate current failures to the preceding revision. It selects the smallest instance $r_G=\operatorname*{arg\,min}_{x\in G}n(x)$ from each failure group, reducing incidental input complexity and making the remaining error easier to isolate. CRS additionally identifies the regression set $\mathcal{U}_t=\mathcal{P}_{t-1}\cap\mathcal{F}_t$, containing instances that were previously feasible but now fail. This comparison also annotates current failures as regressed, persistent, or changed, allowing $d_t$ to emphasize behavior disrupted or left unresolved by the previous revision.

\item \textbf{Scope-Aware Group Representative Sampling (GRS)} concentrates direction evidence on high-priority failure groups. It identifies a focal stage from the stage-level failure distribution and ranks the groups by stage priority and coverage $c_t(G):=|G|/|\mathcal{F}_t|$. It then retains the top $k(\rho_t)$ groups, where fewer groups are retained as feasibility increases, and selects the largest instance from each retained group. Consequently, low-feasibility states expose a broader range of repair needs, whereas high-feasibility states focus the revision on fewer, higher-priority targets.
\end{itemize}

\section{Experiments}
\label{sec:experiments}

In this section, we conduct comprehensive comparison experiments and ablation studies to empirically evaluate the practical effectiveness of \methodname{}.

\begin{table*}[t]
\centering
\caption{ QYI ($\%$) across nine CO benchmark problems. The best, second, and third results in each column are shown in bold, underlined, and marked with $\star$, respectively. }
\vspace{-1em}
\label{tab:problem_results}
\resizebox{0.95\textwidth}{!}{
\begin{tabular}{lcccccccccc}
\toprule
Method & Operator & Intra & Protein & Pedigree & Crew & E-Graph & BPP & CVRP & TSP & Avg. \\
\midrule

EoH & 84.2 & 69.7 & 62.7 & \underline{62.0} & 49.9 & 41.6 & \textbf{99.7} & 85.9 & 49.6 & 67.3 \\

ReEvo & 65.6 & 56.5 & 55.7 & 53.5 & 20.8 & 26.7 & 94.2 & 76.7 & 38.9 & 54.3 \\

HSEvo & 64.7 & 56.9 & 63.9 & 39.9 & 21.9 & 39.1 & \underline{99.5} & 78.0 & 33.9 & 55.3 \\

HiFo-Prompt & 78.6 & 56.8 & 60.5 & \textbf{64.3} & 25.2 & 30.5 & 98.6$^\star$ & 72.4 & 35.3 & 58.0 \\

RefineEvo & 66.3 & 59.2 & 50.5 & 59.8$^\star$ & 35.1 & 41.1 & 87.2 & 74.5 & 35.7 & 56.6 \\

HeuriGym & 91.5$^\star$ & 68.4 & 73.7 & 44.4 & 40.7 & 51.1 & \textbf{99.7} & 85.5 & 52.0 & 67.4 \\

\midrule

DiSH-GWRS & \textbf{93.8} & 69.9$^\star$ & \textbf{81.7} & 57.0 & \underline{67.3} & 56.0$^\star$ & \textbf{99.7} & \underline{93.0} & \underline{85.5} & 78.2$^\star$ \\

DiSH-CRS & \underline{93.6} & \underline{71.2} & \underline{81.5} & 57.7 & 65.5$^\star$ & \textbf{73.2} & \textbf{99.7} & \textbf{93.6} & \textbf{86.9} & \textbf{80.3} \\

DiSH-GRS & 90.5 & \textbf{72.2} & 80.7$^\star$ & 58.0 & \textbf{69.4} & \underline{63.0} & \textbf{99.7} & 92.7$^\star$ & 83.6$^\star$ & \underline{78.9} \\

\bottomrule
\end{tabular}
}
\vspace{-1em}
\end{table*}

\subsection{Experimental Setup}
\label{sec:setup}

We evaluate DiSH-GWRS, DiSH-CRS, and DiSH-GRS on nine problems: six from HeuriGym~\citep{heurigym}, plus bin packing (BPP), capacitated vehicle routing (CVRP), and traveling salesman (TSP) from BPPLIB~\citep{bpplib}, CVRPLIB Set~X~\citep{cvrplibx}, and TSPLIB~\citep{tsplib}, respectively. Demonstration instances guide refinement, while held-out instances are reserved for final evaluation. Baselines include EoH~\citep{eoh}, ReEvo~\citep{reevo}, HSEvo~\citep{hsevo}, HiFo-Prompt~\citep{hifoprompt}, RefineEvo~\citep{refineevo}, and the original HeuriGym refinement procedure. All methods share a complete-solver setting with Qwen3.5-122B-A10B~\citep{qwen35}, DeepSeek-V3~\citep{deepseekv3}, GPT-4.1-mini~\citep{gpt41}, and Gemini-2.5-Flash~\citep{gemini25} as backbones. The five AHD baselines retain their native search operators and procedures. Following HeuriGym, we use the Quality-Yield Index (QYI) as the primary metric and report Quality and Yield separately. Quality uses verified subsets that may differ across methods. Section~\ref{sec:analysis} additionally reports Convergence Rate, Survival Rate, and Tokens to Convergence on demonstration instances. Benchmark details, metric definitions, and baseline adaptations appear in Appendices~\ref{app:benchmarks}, \ref{app:metrics}, and~\ref{app:baseline_configs}, respectively.

\subsection{Comparative Results}
\label{sec:main_results}

\textbf{Per-problem performance.} Across three variants, DiSH achieves the highest QYI on seven of nine problems and ties for best on BPP (Table~\ref{tab:problem_results}). The largest advantages occur on TSP, E-Graph, and Crew Pairing, where the best DiSH variants achieve $86.9$, $73.2$, and $69.4$, respectively. The strongest evaluated baselines on these problems reach $52.0$, $51.1$, and $49.9$. DiSH also leads on Protein Design ($81.7$ vs.\ $73.7$), CVRP ($93.6$ vs.\ $85.9$), Operator Scheduling ($93.8$ vs.\ $91.5$), and Intra ($72.2$ vs.\ $69.7$). On BPP, all three DiSH variants tie with EoH and HeuriGym for the highest QYI of $99.7$. Pedigree remains the exception, with HiFo-Prompt achieving $64.3$ QYI compared with $58.0$ for the best DiSH variant. Appendix~\ref{app:qy_decomposition} further examines Quality and Yield on this benchmark through a complementary comparison with EoH. Overall, the performance advantage extends across most evaluated problems, although the leading DiSH variant varies by task.

\textbf{Overall performance.} Table~\ref{tab:overall} shows that all three DiSH variants outperform every evaluated baseline in mean QYI under each of the four backbones. DiSH-GRS leads on Qwen3.5-122B-A10B ($80.1$), while DiSH-CRS ranks first on DeepSeek-V3 ($80.6$) and GPT-4.1-mini ($82.7$). DiSH-GWRS leads on Gemini-2.5-Flash ($80.9$). DiSH-CRS also achieves the highest overall average QYI of $80.3$, followed by DiSH-GRS ($78.9$) and DiSH-GWRS ($78.2$). Among the baselines, HeuriGym achieves the highest overall average of $67.4$, closely followed by EoH ($67.3$). Overall, the performance advantage extends across all evaluated backbones, although the leading DiSH variant varies with the model. Further analysis in Appendix~\ref{app:qy_decomposition} shows that these QYI gains are accompanied by aggregate improvements in both Quality and Yield for all three DiSH variants. Appendix~\ref{app:problem_backbone} examines how their gains over HeuriGym vary across problems and backbones. Results across shared initializations for all three variants are reported in Appendix~\ref{app:paired}.

\begin{table*}[t]
\centering
\footnotesize
\caption{Overall QYI ($\%$) across backbone models. The best, second, and third results in each column are shown in bold, underlined, and marked with $\star$, respectively. }
\vspace{-1em}
\label{tab:overall}
\begin{adjustbox}{max width=0.85\textwidth}
\begin{tabular}{lccccc}
\toprule
Method & Qwen3.5-122B-A10B & DeepSeek-V3 & GPT-4.1-mini & Gemini-2.5-Flash & Avg. \\
\midrule

EoH & 67.6 $\pm$ 1.5 & 57.3 $\pm$ 6.2 & 68.8 $\pm$ 4.7 & 75.2 $\pm$ 6.8 & 67.3 \\

ReEvo & 46.5 $\pm$ 15.3 & 43.4 $\pm$ 8.3 & 56.8 $\pm$ 8.9 & 70.3 $\pm$ 6.0 & 54.3 \\

HSEvo & 50.0 $\pm$ 11.9 & 44.4 $\pm$ 5.5 & 54.5 $\pm$ 7.1 & 72.4 $\pm$ 2.9 & 55.3 \\

HiFo-Prompt & 60.9 $\pm$ 1.8 & 44.3 $\pm$ 5.4 & 55.6 $\pm$ 3.2 & 71.2 $\pm$ 6.0 & 58.0 \\

RefineEvo & 42.1 $\pm$ 11.2 & 41.6 $\pm$ 4.1 & 67.8 $\pm$ 5.6 & 74.9 $\pm$ 5.7 & 56.6 \\

HeuriGym & 64.1 $\pm$ 5.3 & 62.8 $\pm$ 12.6 & 66.3 $\pm$ 11.2 & 76.6 $\pm$ 9.9 & 67.4 \\

\midrule

DiSH-GWRS & \underline{78.1 $\pm$ 9.0} & 75.9 $\pm$ 8.3$\star$ & 77.9 $\pm$ 7.1$\star$ & \textbf{80.9 $\pm$ 5.7} & 78.2$\star$ \\

DiSH-CRS & 77.6 $\pm$ 7.5$\star$ & \textbf{80.6 $\pm$ 10.4} & \textbf{82.7 $\pm$ 4.9} & \underline{80.4 $\pm$ 7.7} & \textbf{80.3} \\

DiSH-GRS & \textbf{80.1 $\pm$ 7.5} & \underline{76.6 $\pm$ 9.8} & \underline{80.8 $\pm$ 8.6} & 77.9 $\pm$ 10.4$\star$ & \underline{78.9} \\

\bottomrule
\end{tabular}
\end{adjustbox}
\vspace{-1em}
\end{table*}

\subsection{Effects of Direction and Scope}
\label{sec:analysis}

\begin{wraptable}{r}{0.50\linewidth}
    \centering
    \small
    \vspace{-2em}
    \caption{Ablation of direction and scope control.}
    \label{tab:ablation}
    \vspace{-1em}
    \begin{adjustbox}{max width=0.95\linewidth}
    \begin{tabular}{@{}lcccc@{}}
        \toprule
        Method & \shortstack{QYI \\(\%) $\uparrow$} & \shortstack{Conv.\\(\%) $\uparrow$} & \shortstack{Survival\\(\%) $\uparrow$} & \shortstack{Tokens to\\Conv. $\downarrow$} \\
        \midrule
        HeuriGym & 67.4 & 55.2 & 40.0 & 19{,}897 \\
        DiSH (w/o Sco) & 69.0 & 55.6 & 46.6 & 16{,}335 \\
        DiSH (w/o Dir) & 76.4 & \underline{91.4} & 92.3$^\star$ & 14{,}804 \\
        DiSH-GWRS & 78.2$^\star$ & 88.3 & \underline{93.0} & 11{,}270$^\star$ \\
        \midrule
        DiSH-CRS & \textbf{80.3} & \textbf{92.9} & \textbf{94.2} & \textbf{10{,}271} \\
        DiSH-GRS & \underline{78.9} & 91.4$^\star$ & 91.9 & \underline{10{,}868} \\
        \bottomrule
    \end{tabular}
    \end{adjustbox}
    \vspace{-1em}
\end{wraptable}
Table~\ref{tab:ablation} shows larger gains in convergence and survival from scope control, while combining scope and direction further improves QYI and reduces Tokens to Convergence. We conduct a controlled $2 \times 2$ ablation with all other settings fixed. HeuriGym enables neither control; DiSH (w/o Sco) uses GWRS feedback alone; DiSH (w/o Dir) retains full feedback and enables scope control; and DiSH-GWRS combines both. DiSH-CRS and DiSH-GRS additionally test alternative direction strategies with both controls enabled. Evaluation includes held-out test QYI and three diagnostics on demonstration instances, defined in Appendix~\ref{app:trajectory_metrics}. Compared with HeuriGym, DiSH (w/o Dir) raises Convergence Rate from 55.2\% to 91.4\% and Survival Rate from 40.0\% to 92.3\%. DiSH (w/o Sco) yields smaller increases to 55.6\% and 46.6\%, respectively. Both configurations improve QYI and reduce Tokens to Convergence relative to HeuriGym. All three combined variants further improve these two metrics over either control alone, although convergence and survival do not improve uniformly. Relative to (DiSH (w/o Dir)), DiSH-GWRS has a lower Convergence Rate ($88.3\%$ versus $91.4\%$), and DiSH-GRS has a slightly lower Survival Rate ($91.9\%$ versus $92.3\%$). Overall, joint control improves final performance and token efficiency, with DiSH-CRS achieving the best results across all four metrics.

\begin{figure*}[t]
    \centering
    \vspace{-1em}
    \includegraphics[width=1.0\textwidth]{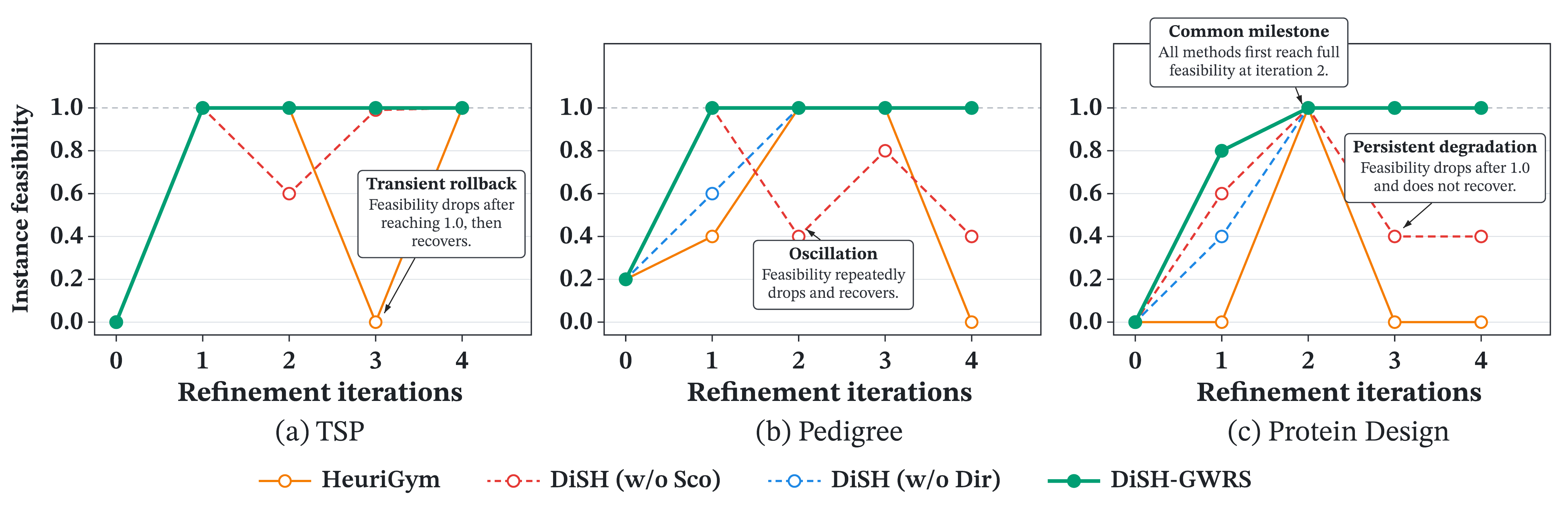}
    \vspace{-2em}
    \caption{Representative feasibility trajectories on TSP, Pedigree, and Protein Design. }
    \vspace{-1em}
    \label{fig:traj}
\end{figure*}

\textbf{Feasibility stability across refinement.} Representative trajectories in Figure~\ref{fig:traj} illustrate the stabilizing effect of adaptive scope control. On TSP, HeuriGym and DiSH (w/o Sco) temporarily lose feasibility before recovering. On Pedigree, DiSH (w/o Sco) oscillates, while HeuriGym drops from full feasibility to zero in the final iteration. On Protein Design, all methods first reach full feasibility at iteration $2$, but HeuriGym and DiSH (w/o Sco) subsequently regress without recovery. In contrast, DiSH (w/o Dir) and DiSH-GWRS maintain full feasibility once reached in all three examples, with their curves overlapping on TSP. These patterns align with their higher Survival Rates in Table~\ref{tab:ablation}, suggesting that adaptive scope control based on feasibility helps reduce regressions.

\subsection{Controlled Mechanism Analysis}
\label{sec:controlled_analysis}

\begin{wrapfigure}{r}{0.4\linewidth}
    \centering
    \vspace{-1.8em}
    \includegraphics[width=\linewidth]{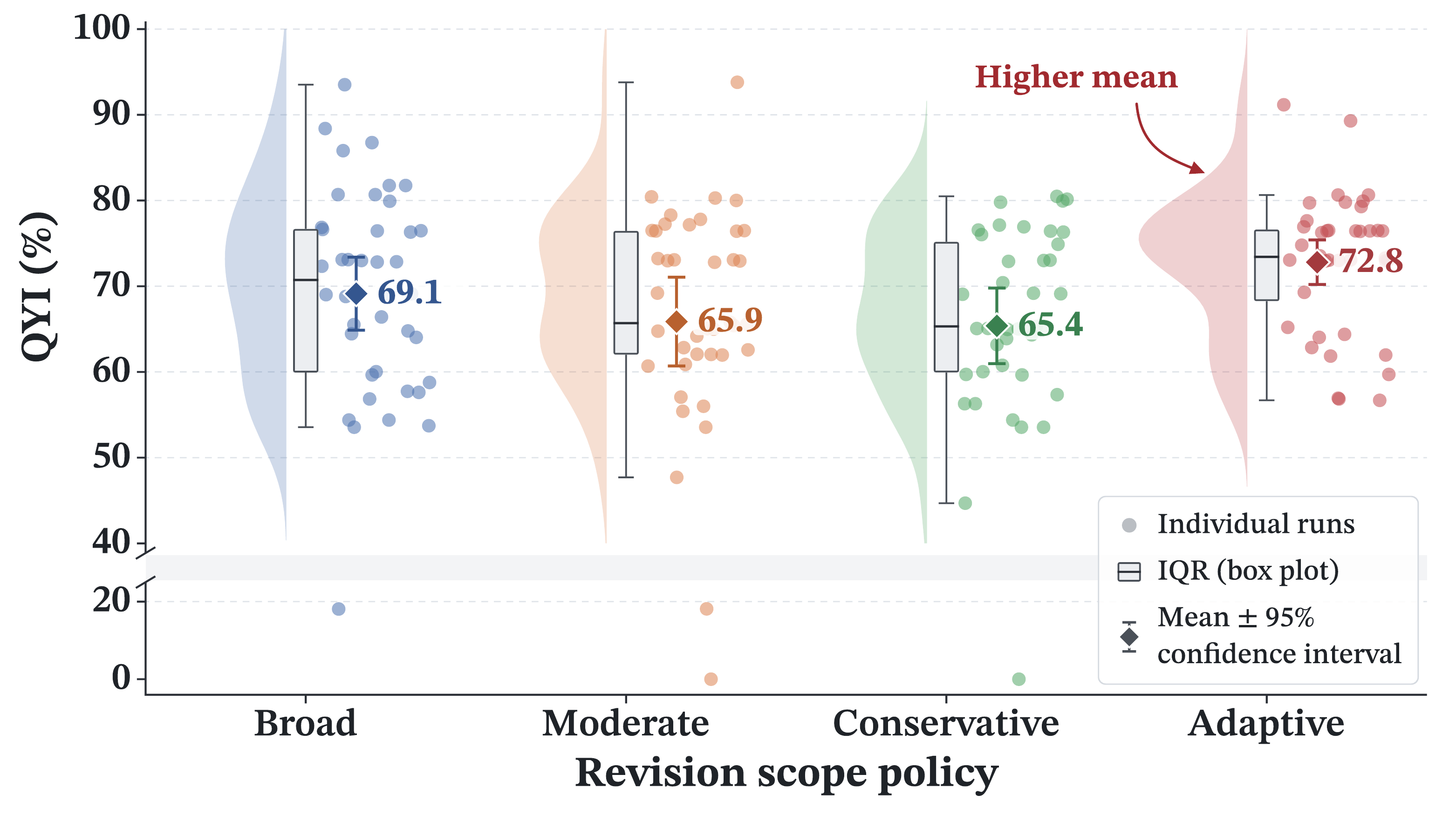}
    \vspace{-2em}
    \caption{Adaptive vs.\ fixed scopes.}
    \label{fig:adaptive_fixed_revision_scopes}
    \vspace{-1.5em}
\end{wrapfigure}
\textbf{Adaptive versus fixed revision scopes.} We compare adaptive and fixed revision scopes to test whether the observed gains can be explained by a single fixed revision scope. The only difference is the scope policy: Adaptive selects the revision scope according to the current feasibility state, whereas the fixed baselines use Broad, Moderate, or Conservative scope throughout refinement. Figure~\ref{fig:adaptive_fixed_revision_scopes} reports the QYI distribution and mean for each policy. Adaptive achieves the highest average QYI ($72.8\%$), outperforming Broad ($69.1\%$), Moderate ($65.9\%$), and Conservative ($65.4\%$), while also exhibiting lower variability than the fixed-scope baselines. These results suggest that the benefit comes not from any single fixed scope, but from adapting revision scope to the current feasibility state.

\begin{wrapfigure}{r}{0.4\linewidth}
    \centering
    \vspace{-0.8em}
    \includegraphics[width=\linewidth]{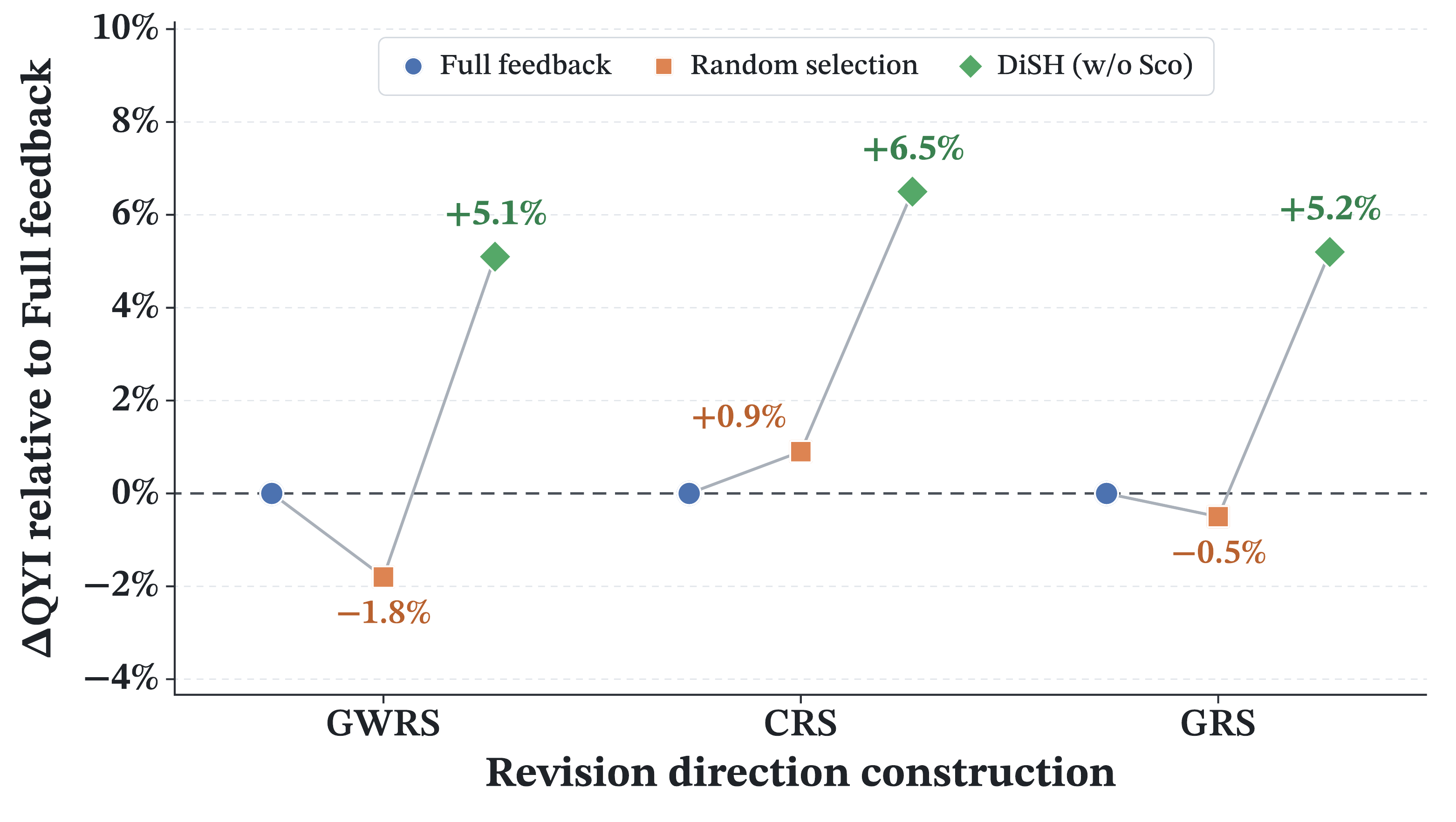}
    \vspace{-2em}
    \caption{Evidence organization.}
    \label{fig:structured_revision_direction}
    \vspace{-1.5em}
\end{wrapfigure}
\textbf{Structured revision direction.} We compare structured direction construction with full feedback and random selection to distinguish evidence organization from feedback reduction alone. Random selection uses the same evidence budget as the corresponding structured strategy. Adaptive scope control is disabled, while revision scope and all other settings remain fixed. Figure~\ref{fig:structured_revision_direction} reports QYI changes relative to full feedback. Structured direction construction achieves gains of $5.1\%$, $6.5\%$, and $5.2\%$ for GWRS, CRS, and GRS, respectively. Random selection yields changes of $-1.8\%$, $+0.9\%$, and $-0.5\%$, showing no consistent improvement. All three strategies outperform random selection with the same evidence budget. This suggests that the gains come from organizing failure evidence, not just reducing feedback.

\begin{figure*}[t]
    \vspace{-1em}
    \centering
    \vspace{-1em}
    \includegraphics[width=1.0\textwidth]{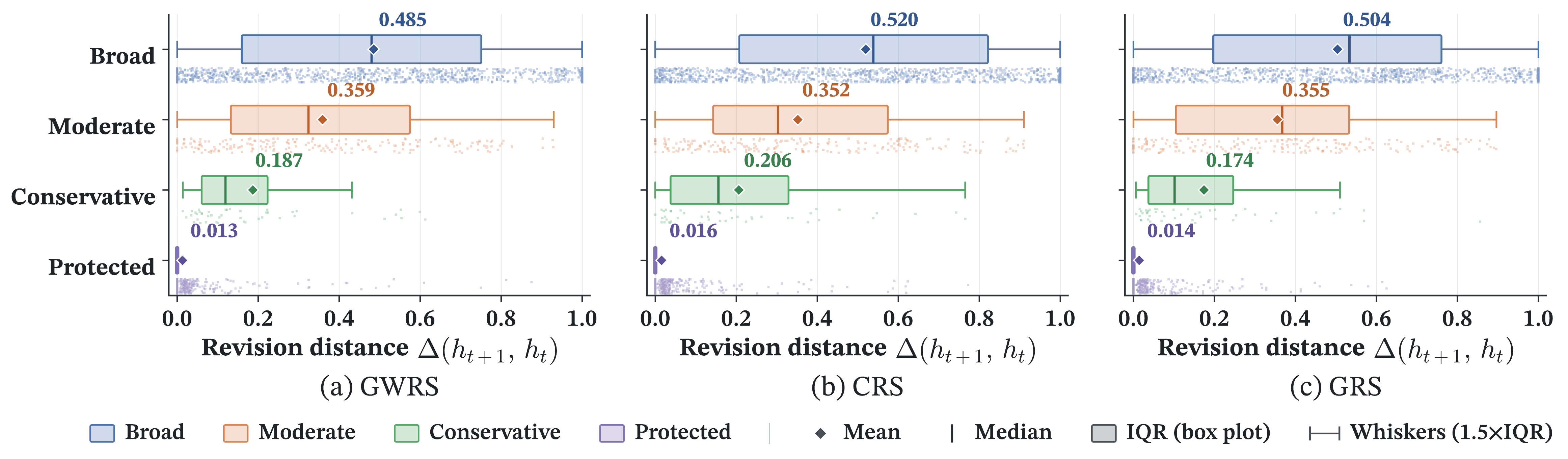}
    \vspace{-2em}
    \caption{Revision distance by semantic scope.}
    \label{fig:scope_distance}
\end{figure*}

\textbf{Joint Effects of Revision Direction and Scope.} Across all three direction strategies, narrower semantic scopes are associated with smaller mean revision distances. Figure~\ref{fig:scope_distance} shows a consistent decrease from \textsc{Broad} through \textsc{Moderate} and \textsc{Conservative} to \textsc{Protected}. Mean distances range from $0.485$ to $0.520$ under \textsc{Broad} and approach zero under \textsc{Protected} ($0.013$ to $0.016$). These results indicate that scope guidance before generation is reflected in progressively smaller code changes across different direction constructions. This consistent pattern across direction strategies provides empirical support for the scope-control mechanism described in Section~\ref{subsec:scope}.

\subsection{Case Studies}
\label{sec:scope_cases}

\begin{figure*}[t]
    \centering
    \vspace{-1em}
    \includegraphics[width=1.0\textwidth]{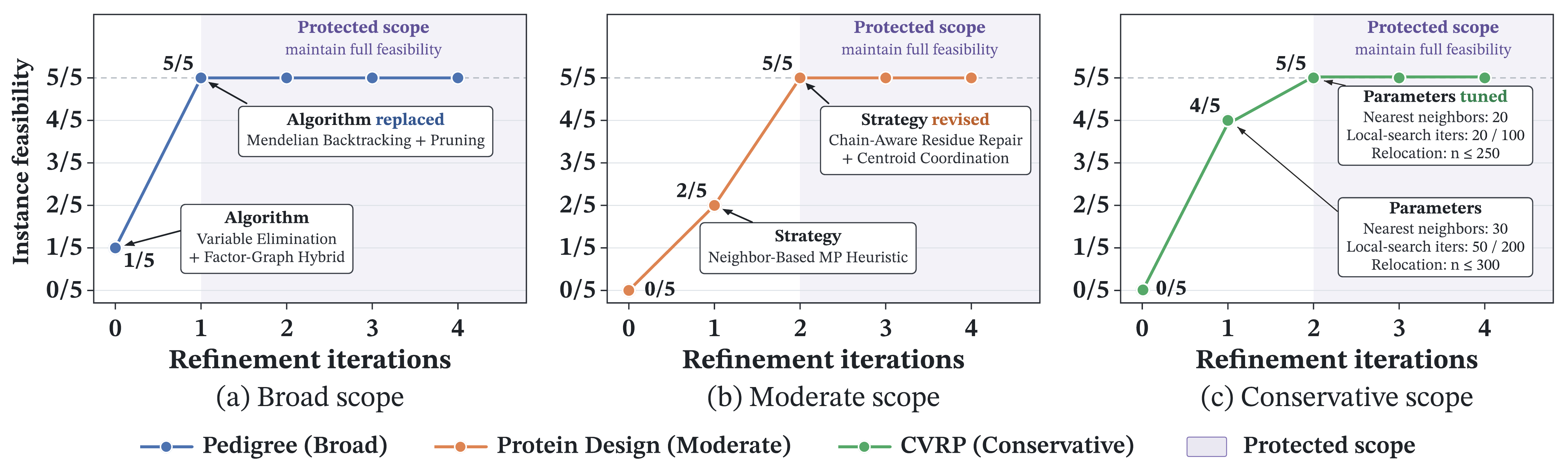}
    \vspace{-2em}
    \caption{ Representative solver revisions under different semantic revision scopes.}
    \vspace{-0.5em}
    \label{fig:scope_case}
\end{figure*}

Figure~\ref{fig:scope_case} presents three representative revisions from DiSH runs, illustrating how different scope modes guide solver modification. On Pedigree, only one of five demonstration instances passes, so \textsc{Broad} mode permits aggressive changes to the core algorithm. The revision replaces the variable elimination and factor graph hybrid with Mendelian backtracking and pruning, raising feasibility from $1/5$ to $5/5$. On Protein Design, \textsc{Moderate} mode applies at $2/5$ feasibility. The revision introduces chain-aware residue repair and centroid coordination, raising demonstration feasibility from $2/5$ to $5/5$. On CVRP, four of five instances already pass, and the \textsc{Conservative} revision makes targeted parameter adjustments. It reduces the neighbor count, local search iterations, and relocation threshold without changing the routing or search logic, raising feasibility from $4/5$ to $5/5$. The purple regions mark \textsc{Protected} mode, entered after full feasibility is reached. This mode preserves the solver by default and permits only parameter changes with a concrete benefit and no adverse effect on currently feasible instances. All three runs maintain full feasibility throughout the remaining iterations. Together, these cases illustrate how DiSH adapts revision scope to feasibility, allowing substantial repairs when needed and prioritizing preservation as feasibility improves.

\section{Conclusion}

This work shows that stable LLM-based refinement of combinatorial optimization solvers depends not only on the feedback itself, but also on how that feedback is translated into solver revisions. Existing approaches largely leave two key decisions to the LLM: \emph{what to change} and \emph{how much to change}. This can cause revisions to disrupt previously feasible behavior and lead to feasibility oscillation. \methodname{} makes these decisions explicit by controlling revision direction and scope: representative failure groups identify what should be addressed, while the current feasibility state determines how extensively the solver should be revised. Across nine combinatorial optimization problems and four backbone LLMs, \methodname{} achieves the highest overall average QYI, produces more stable refinement trajectories, and uses roughly half the feedback tokens required by full feedback. These results suggest that controlling the revision process itself provides another effective way to improve automatic heuristic design, complementing advances in search and solver generation.

\section*{AI Use Statement}

In this work, we used generative AI tools solely for language polishing, grammar checking, and discussions of experimental coverage. All AI-assisted content was reviewed by the authors, who take full responsibility for the final content of this work.

\section*{Ethics Statement}

This work studies LLM-based solver refinement on combinatorial optimization benchmarks. The experiments focus on benchmark evaluation, and practical deployment would require further validation under application-specific requirements.

\section*{Reproducibility Statement}

The method and refinement procedure are described in Section~\ref{sec:method} and Appendix~\ref{app:harness-loop}. Appendix~\ref{app:prompts} provides the prompts and refinement examples. Metric definitions and experimental configurations are detailed in Appendices~\ref{app:metrics} and~\ref{app:experimental_setup}, with additional results reported in Appendix~\ref{app:additional_results}.
\bibliographystyle{iclr2027_conference}
\bibliography{references}

\clearpage
\appendix

\startcontents[appendix]

\thispagestyle{plain}

\begin{center}
    {\Large\bfseries Appendix Contents}
\end{center}

\vspace{1em}

\setcounter{tocdepth}{2}
\printcontents[appendix]{}{1}{}

\clearpage

\section{Related Work}
\label{sec:related}

\subsection{Harness for Language-Model Agents}
\label{sec:related-harness}
Before the agent harness emerged as an explicit systems concept, external environments already supported agent interaction and feedback. TextWorld~\citep{cote2019textworld}, ALFWorld~\citep{alfworld}, and WebShop~\citep{webshop} provide executable environments for text, embodied, and web tasks. ReAct~\citep{react} integrates reasoning, action, and observation, while Reflexion~\citep{reflexion} carries failure feedback across trials through verbal reflection. Research has expanded to execution, verification, and interface design. SWE-bench~\citep{swebench} tests repository edits, while WebArena~\citep{webarena} checks functional correctness against environment states. SWE-agent~\citep{sweagent} shows that interface design can improve performance without changing the model. Together, these developments mark a transition from interaction environments to unified harnesses integrating execution, tooling, context, state maintenance, and verification.

Recent work further uses execution histories and failures to organize subsequent runs. Life-Harness~\citep{lifeharness} distills reusable interventions from recurring failures in training interaction histories while keeping the harness fixed during evaluation. AgentTether~\citep{agenttether} localizes critical subtrajectories within failed executions and provides targeted guidance and runtime interventions for re-execution of the same agent. These works illustrate a shift from using execution outcomes for evaluation toward using them to organize subsequent execution.

\subsection{LLMs for Optimization Modeling and Solving}
\label{sec:related_modeling}

LLMs are increasingly used to automate optimization modeling by translating natural-language problem descriptions into mathematical formulations and solver-executable code. Training-based approaches develop models specialized for this task through synthesized supervision. ORLM combines semi-automated data generation with open-source model training~\citep{orlm}, while OptMATH constructs scalable bidirectional modeling data with forward verification and rejection sampling~\citep{optmath}.SIRL further grounds training in solver-verifiable signals over generated models and executable code, providing rewards for syntactic validity, feasibility, and solution quality~\citep{sirl}. Opt-Miner combines tree-guided data synthesis with reinforcement learning that rewards retrieval quality and modeling accuracy, training agents to acquire external knowledge for optimization problems with knowledge gaps~\citep{liu2026optminer}.

Complementary prompt-based and agentic approaches structure the modeling process itself. OptiMUS formulates models, generates solver code, and uses execution feedback to diagnose and repair modeling or implementation errors~\citep{ahmaditeshnizi2024optimus}. OptiTree organizes formulation through a hierarchical modeling tree~\citep{optitree}, while MM-Agent studies broader mathematical modeling through staged workflows and hierarchical knowledge~\citep{mmagent}. LEAN-LLM-OPT dynamically constructs agentic workflows for large-scale auto-formulation and delegates structured data operations to auxiliary tools~\citep{leanllmopt}. MIRROR combines hierarchical exemplar retrieval with execution-driven revision of mathematical models and solver code~\citep{mirror}. Evaluation has also advanced through structural-equivalence measures~\citep{orgeval} and industrial-scale benchmarks derived from real mixed-integer programs~\citep{miplibnl}.

Beyond formulation, LLMs have been trained as end-to-end CO solvers~\citep{llmcosolver} and used to automate solver-level decisions, including instance-specific MIP configuration~\citep{grimip} and branching-policy discovery~\citep{llm4branch}. These studies focus on problem formulation, direct solution generation, or solver-level decisions. A related line of work instead uses LLMs to automatically design and refine executable heuristics.

\subsection{LLM-based Automatic Heuristic Design}
\label{sec:related_ahd}

Heuristic design for combinatorial optimization has long relied on expert knowledge. Hyper-heuristics and algorithm configuration methods search within predefined spaces of heuristics, components, or parameters~\citep{burke2013,pillay2018,burke2019,paramils,irace}. LLMs extend this search to executable programs~\citep{liu2026survey}. FunSearch~\citep{funsearch} evolves functions within fixed program skeletons, EoH~\citep{eoh} evolves heuristic ideas and implementations, and ReEvo~\citep{reevo} incorporates reflection. Together, these methods establish \emph{iterative revision} through repeated program evaluation and verification.

Subsequent research focuses on organizing and regulating this iterative search process. Approaches promote candidate diversity and branching~\citep{hsevo,mctsahd}, introduce different modification granularities~\citep{a2dept,tide}, and adapt search behavior to current states~\citep{hifoprompt,refineevo,esahd}. Complementary extensions explore complete metaheuristics and coupled operators~\citep{llamea,glns}, collaborative search~\citep{partevo}, multiple objectives~\citep{meoh}, and complementary heuristic sets~\citep{eohs}. Others improve guidance and evaluation efficiency~\citep{hercules} or adapt the generating model~\citep{calm}. These developments broaden both the algorithmic structures explored and the mechanisms governing their evolution.

Alongside these advances in search design, other work focuses on enriching and organizing the evidence used to guide revisions. HeuriGym~\citep{heurigym}, LaGO~\citep{lago}, and ReVEL~\citep{revel} provide execution, analytical, and structured feedback, respectively. PathWise~\citep{pathwise} reuses search histories, AHD Agent~\citep{ahdagent} acquires targeted evidence, and ARES~\citep{ares} combines structural ablation with accumulated design experience. These approaches extend revision guidance beyond aggregate performance toward richer, more structured, and more targeted evidence.

\subsection{Harnesses in LLM-based Solver Refinement}
\label{sec:related_harness}
For solver refinement, harness design must address not only evaluation and feedback, but also how each revision is controlled. HeuriGym~\citep{heurigym} integrates execution, output parsing, and verification to provide instance-level feedback. CO-Bench~\citep{cobench} and ALE-Bench~\citep{alebench} evaluate complete solvers through execution and performance measurement. However, mapping evaluation states to explicit revision controls remains underexplored.

This gap concerns both \emph{revision scope} and \emph{revision direction}. Scope determines how extensively a solver may change. Existing methods typically couple it to search operators rather than constrain it independently before generation~\citep{mctsahd,a2dept,hifoprompt,refineevo}. Direction specifies where changes should focus, yet priorities largely depend on the model's interpretation of richer feedback~\citep{lago,revel,pathwise,ahdagent}. The practical challenge is to repair current failures without disrupting previously verified behavior. This motivates exploring how a harness can use the execution and verification evidence it collects to jointly control revision direction and scope.

\section{Method Details}\label{app:harness-loop}
\subsection{Refinement Procedure}
DiSH follows a harness-controlled iterative refinement procedure, summarized in Algorithm~\ref{alg:dish_refinement}. The process begins by evaluating the initial solver $h_0$ on the demonstration set to obtain the evaluation records $\mathcal{Z}_0$ and the corresponding feasibility pass rate $\rho_0$. At each refinement step $t$, the harness first determines the revision scope $\eta_t$ from the current feasibility pass rate $\rho_t$ and constructs the revision direction $d_t$ from the accumulated evaluation records $\mathcal{Z}_{0:t}$. These controls, together with the evaluation feedback, are assembled into the refinement prompt $p_t$. The complete LLM context is then defined as $\mathcal{C}_t=[\mathcal{S},\mathcal{U},h_t,p_t]$, consisting of the fixed system prompt, task prompt, current solver, and round-specific refinement prompt. The LLM uses this context to generate the next solver $h_{t+1}$, which is evaluated in the same manner to obtain $\mathcal{Z}_{t+1}$ and $\rho_{t+1}$. All experiments use \texttt{history\_rounds=1}, so each LLM call retains only this controlled context rather than the full interaction history. After $K$ refinement steps, all evaluated solvers remain eligible for delivery, and DiSH selects as $h^\star$ the checkpoint with the highest demonstration feasibility pass rate $\rho_t$, breaking ties by lower mean cost over the feasible demonstration instances. The selected solver is then fixed and evaluated on the held-out test instances, which are not used during refinement or checkpoint selection.

\begin{algorithm}[htbp]
\caption{Direction-and-Scope-Guided Solver Refinement}
\label{alg:dish_refinement}
\begin{algorithmic}[1]
\Require system prompt $\mathcal{S}$; task prompt $\mathcal{U}$ (Appendix~\ref{app:prompts}); demonstration set $\mathcal{D}_{\mathrm{train}}$; initial solver $h_0$; direction strategy $\mathcal{G}_{\mathrm{dir}} \in \{\mathcal{G}_{\mathrm{GWRS}}, \mathcal{G}_{\mathrm{CRS}}, \mathcal{G}_{\mathrm{GRS}}\}$; number of refinement steps $K$
\Ensure delivered solver $h^\star$

\State $(\mathcal{Z}_0,\rho_0) \gets \mathcal{E}(h_0,\mathcal{D}_{\mathrm{train}})$

\For{$t=0,\ldots,K-1$}
    \State $\eta_t \gets \mathcal{G}_{\mathrm{step}}(\rho_t)$
    \State $d_t \gets \mathcal{G}_{\mathrm{dir}}(\mathcal{Z}_{0:t})$
    \State $p_t \gets \mathcal{H}(\mathcal{Z}_t,d_t,\eta_t)$
    \State $\mathcal{C}_t \gets [\mathcal{S},\mathcal{U},h_t,p_t]$
    \State $h_{t+1} \sim \mathcal{M}_{\mathrm{LLM}} (\cdot\mid\mathcal{C}_t)$
    \State $(\mathcal{Z}_{t+1},\rho_{t+1}) \gets \mathcal{E}(h_{t+1},\mathcal{D}_{\mathrm{train}})$
\EndFor

\State $h^\star \gets \operatorname{Deliver} \left( \{(h_t,\mathcal{Z}_t)\}_{t=0}^{K} \right)$
\State \Return $h^\star$

\end{algorithmic}
\end{algorithm}

\subsection{Evaluation and Scope Rules}\label{app:evaluation_scope}
\subsubsection{Instance Evaluation}
\label{app:instance_evaluation}

The evaluation module $\mathcal{E}$ assesses each solver on the demonstration instances through three sequential stages. An instance proceeds to the next stage only after passing the current one. If a stage fails, the evaluation record identifies this first failure stage and the corresponding error. Instances passing all three stages are marked as \texttt{SUCCESS}. Table~\ref{tab:stages} summarizes the evaluation criteria and failure outcomes. These records distinguish execution failures, invalid outputs, and constraint violations. They are aggregated into the harness state diagnosis described in Appendix~\ref{app:diagnosis}, which summarizes the current solver state for subsequent revision.

\begin{table}[h]
\centering
\small
\caption{Three-stage feasibility evaluation for each instance.}
\label{tab:stages}
{
\begin{tabular}{@{}l p{10.2cm}@{}}
\toprule
Stage & Evaluation criterion \\
\midrule
I\ \ Execution & Terminates normally within the time and memory limits; otherwise returns an exception, timeout, or memory error. \\

II\ \ Output validity & Produces a non-empty output and allows the evaluation pipeline to complete; otherwise reports an output or evaluation error. \\

III\ \ Solution feasibility & Receives a positive verdict from the problem-specific verifier; otherwise reports a verification failure. \\
\bottomrule
\end{tabular}
}
\end{table}

\subsubsection{Revision Scope Rules}
\label{app:scope_rules}

The scope controller maps the current feasibility pass rate $\rho_t$ to a revision instruction $\eta_t$. Table~\ref{tab:scope_modes} specifies the four modes, their thresholds, and the permitted revision levels. The corresponding prompts are provided in Appendix~\ref{app:stepsize}.

\begin{table}[t]
\centering
\caption{Feasibility-adaptive revision-scope instructions.}
\label{tab:scope_modes}
\small
\begin{tabular}{@{}l l p{9cm}@{}}
\toprule
Feasibility state & Instruction & Permitted revision scope \\
\midrule

$\rho_t < 0.4$ & \textsc{Broad} & Parameter-, strategy-, and algorithm-level changes are permitted, including algorithm replacement. \\

$0.4 \leq \rho_t < 0.8$ & \textsc{Moderate} & Strategy-level changes are preferred; algorithm-level changes are permitted only when the observed failures indicate a structural mismatch. \\

$0.8 \leq \rho_t < 1$ & \textsc{Conservative} & Parameter-level changes and narrowly targeted strategy-level changes are permitted; algorithm-level changes are prohibited. \\

$\rho_t = 1$ & \textsc{Protected} & The current solver is preserved by default; only parameter-level changes with a concrete benefit and no adverse effect on currently feasible instances are permitted. \\

\bottomrule
\end{tabular}
\end{table}

\section{Prompt Design}\label{app:prompts}
In this section, we detail the system and user prompts used by DiSH, together with the harness state diagnosis and the controls over revision scope and direction. We also provide a complete refinement example to illustrate how these components work together in practice. In the prompt text, ``test case'' is the literal label inherited from the HeuriGym feedback template and refers exclusively to demonstration instances in $\mathcal{D}_{\mathrm{train}}$; held-out test instances are never exposed to the LLM.

\subsection{System Prompts}\label{app:sysprompt}
Each refinement iteration begins with a task-agnostic system prompt that instructs the LLM to revise the current solver under a fixed set of execution, output, and refinement constraints. Beyond the basic solver interface and runtime requirements, the prompt specifies three revision levels, the three-stage harness evaluation pipeline, and a feasibility-conditioned scope schedule that progressively restricts how much the solver may change as more demonstration cases succeed. Several placeholders are instantiated at runtime according to the problem configuration. For example, \ph{TIMEOUT} specifies the wall-clock time limit, and \ph{N\_DEMO} denotes the number of demonstration instances used for evaluation feedback.

\begin{promptbox}{System Prompt}{sysprompt}
\small

You are an expert solver engineer. Iteratively refine a correct, robust solution to the optimization problem below: keep what works and repair failures in a targeted way. How large a change is appropriate depends on how much of the solver already works---this is governed by the mode schedule in point 7, not by a blanket preference for small changes.

\medskip
Follow these instructions strictly.

\begin{promptlist}

\item[1.] Implement only the \texttt{solve} function. Do NOT change its signature, argument names, or argument types. Helper functions, classes, and internal data structures are allowed.

\item[2.] The machine has 8 CPU cores. Time limit: \ph{TIMEOUT} seconds. Use parallelism only when it is reliable and clearly beneficial within the budget.

\item[3.] Python 3.12. Permitted third-party libraries only:
  \begin{promptlist}
  \item[--] \texttt{numpy==2.2.5}
  \item[--] \texttt{networkx==3.4.2}
  \item[--] \texttt{pandas==2.2.3}
  \end{promptlist}

\item[4.] Respond with code only---a complete \texttt{solve} implementation plus any helpers. No explanation, no Markdown, no extra text.

\item[5.] \textbf{Change scope}---three levels, distinguished by how much of the solver a change touches, not by which component is edited:
  \begin{promptlist}

  \item[--] \textbf{Parameter}: numerical constants only; control flow and logic unchanged.

  \item[--] \textbf{Strategy}: logic changes within the current algorithmic paradigm; approach and complexity class unchanged.

  \item[--] \textbf{Algorithm}: replaces the core paradigm or changes the complexity class.

  \end{promptlist}

Correctness fixes---a wrong loop bound, a wrong index, a mishandled constraint, an off-by-one---are NOT one of these three levels and are never restricted by mode. Fix them in any state, including when the defect lies in code that currently passing cases also execute.

\item[6.] Each round you receive concrete execution feedback computed from a fixed set of \ph{N\_DEMO} demonstration cases: per-case results, the dominant failure stage, and the demo SUCCESS proportion. The stages are:
  \begin{promptlist}

  \item[--] \textbf{Stage I}: execution failure (crash, exception, import error, or timeout);

  \item[--] \textbf{Stage II}: the program runs but produces missing or malformed output;

  \item[--] \textbf{Stage III}: output is produced but rejected by the verifier (infeasible).

  \end{promptlist}

This feedback is your primary signal---act on it; do not refine from intuition alone. The \ph{N\_DEMO} cases are representative but not exhaustive, and the final solver is judged on unseen instances, so prefer general robustness over fitting these \ph{N\_DEMO}. Breaking a currently passing case is a regression.

\item[7.] The permitted change scope narrows as demo SUCCESS rises---broad when little works, tight when most works. You will be told at the start of each round which mode is active and exactly what it permits; obey that instruction. This clause only states the overall shape:
  \begin{promptlist}

  \item[--] \textbf{$<40\%$ (BROAD)}: all three levels open. Restore the dominant failure stage first and do not optimize cost yet. Even if a few cases already pass, treat this as full recovery of that stage---those cases may simply not exercise the faulty path.

  \item[--] \textbf{$40$--$80\%$ (MODERATE)}: Strategy by default; escalate to Algorithm only when a failing case points to a structural mismatch in the current paradigm; Parameter anytime.

  \item[--] \textbf{$80$--$99\%$ (CONSERVATIVE)}: Parameter; Strategy only when narrowly targeted at the specific failing cases; Algorithm not permitted.

  \item[--] \textbf{$100\%$ (PROTECTED)}: effectively frozen. Return the solver unchanged unless a Parameter change is both provably safe on every passing case and clearly beneficial.

  \end{promptlist}

Global rules (apply in every mode):
  \begin{promptlist}

  \item[--] Address the failure stage with the most failing cases first (ties go to the earlier stage). Overall priority: executable $\rightarrow$ valid output $\rightarrow$ feasible $\rightarrow$ lower cost. Never sacrifice an earlier stage to chase a later one.

  \item[--] After any change, every previously passing case must still pass---judged by the end result, not by treating any block of code as untouchable.

  \item[--] Correctness fixes (point 5) are permitted in every mode without restriction.

  \end{promptlist}

\item[8.] Respect the exact input/output format in the problem statement, handle edge cases, and keep input parsing and output writing stable.

\item[9.] Do not hardcode behavior for individual demonstration cases; the solver must generalize to unseen inputs.

\end{promptlist}

\medskip
Above all: keep each change within the scope the active mode permits---broad when little works, tight when most works. At 100\%, usually change nothing. Correctness fixes are always allowed.

\end{promptbox}

\subsection{User Prompts}\label{app:userprompt}
For each problem, the first refinement iteration begins with a task-specific user prompt whose structure is adapted from the HeuriGym setup~\citep{heurigym}. The prompt provides the problem information, which specifies the optimization task together with its background, formalization, and input and output requirements, as well as a program template that defines the solver interface used throughout refinement. Accordingly, \ph{BACKGROUND}, \ph{FORMALIZATION}, \ph{INPUT\_FORMAT}, and \ph{OUTPUT\_FORMAT} are instantiated from the corresponding problem description, while \ph{SOLVER\_TEMPLATE} is filled with the problem-specific solver template.

\begin{promptbox}{User Prompt}{userprompt}
\small

\texttt{\# Problem Information}

\medskip
\texttt{\#\# Background}

\ph{BACKGROUND}

\medskip
\texttt{\#\# Formalization}

\ph{FORMALIZATION}

\medskip
\texttt{\#\# Input Format}

\ph{INPUT\_FORMAT}

\medskip
\texttt{\#\# Output Format}

\ph{OUTPUT\_FORMAT}

\medskip
\texttt{\# Program Template}

\ph{SOLVER\_TEMPLATE}

\end{promptbox}

\subsection{Harness State Diagnosis}\label{app:diagnosis}

After each Harness evaluation, the instance-level evaluation records $Z_t$ of the current solver are organized into a round-level \emph{state diagnosis}. For every evaluated solver state, the diagnosis reports the numbers and proportions of Stage~I, Stage~II, Stage~III, and SUCCESS cases, together with the current \texttt{Focus stage}, \texttt{Current mode}, and the corresponding \texttt{Mode reason}. When sufficient refinement history is available, the refinement prompt additionally includes a compact \texttt{Progress Trend} summarizing the preceding evaluated solver states. These quantities help the LLM identify the remaining failure patterns, assess how much solver behavior has already been verified, and track changes across refinement rounds. The same structured state also provides the basis for the subsequent revision-scope control and revision-direction construction. The resulting diagnosis is placed at the beginning of each round-specific refinement prompt, providing the LLM with a concise view of the current solver state and its recent refinement progress.

\begin{promptbox}{Harness State Diagnosis}{diagnosis}
\small\ttfamily\raggedright

\# Feedback from Previous Iteration (Iteration \ph{ITERATION})

\medskip
\#\# Demo/Train Results Summary

\medskip
- Total demo/train cases: \ph{N\_DEMO}\\
- Stage I failures (execution error): \ph{N\_I} (\ph{R\_I})\\
- Stage II failures (no/invalid output): \ph{N\_II} (\ph{R\_II})\\
- Stage III failures (verifier rejected): \ph{N\_III} (\ph{R\_III})\\
- Successful cases: \ph{N\_SUCCESS} (\ph{R\_SUCCESS})\\
- Focus stage: \ph{FOCUS\_STAGE}\\
- Current mode: \ph{MODE}\\
- Mode reason: \ph{MODE\_REASON}

\medskip
\textit{When sufficient refinement history is available:}

\medskip
\#\# Progress Trend (demo/train)

\medskip
\ph{PROGRESS\_TREND}

\end{promptbox}

Specifically, SUCCESS proportion corresponds to the feasibility pass rate $\rho_t$ defined in the main text. The revision scope then maps $\rho_t$ to \texttt{Current mode}. The \texttt{Focus stage} identifies the evaluation stage with the largest number of failed cases $f_t=\operatorname*{arg\,max}_{s\in\{\mathrm{I},\mathrm{II},\mathrm{III}\}}\left|\left\{x\in D_{\mathrm{train}}:s_t(x)=s\right\}\right|$. Ties are resolved by prioritizing earlier stages in the evaluation pipeline: Stage~I $\rightarrow$ Stage~II $\rightarrow$ Stage~III. When all demonstration cases succeed, no failure stage is active and the \texttt{Focus stage} is reported as \texttt{none (all demonstration cases pass)}.

\subsection{Revision Scope Prompts}\label{app:stepsize}
Revision scope control determines how extensively the LLM is allowed to modify the current solver at each refinement round. Early in refinement, substantial changes may be necessary to repair fundamental failures, whereas later revisions should increasingly preserve behavior that has already been verified. DiSH therefore adapts the permitted revision scope to the current solver state by mapping the SUCCESS proportion $\rho_t$ to one of four modes, \textsc{Broad}, \textsc{Moderate}, \textsc{Conservative}, or \textsc{Protected}. As $\rho_t$ increases, these modes progressively narrow the allowable changes from unrestricted Parameter, Strategy, and Algorithm level revisions to increasingly localized modifications and, ultimately, preservation of the current solver.

\textbf{\textsc{Broad}.} When $\rho_t<0.4$, the solver remains in an early recovery stage, where failures may arise from execution, output generation, or solution feasibility. DiSH therefore permits Parameter, Strategy, and Algorithm-level revisions without restricting the change scope. The guidance is further specialized by the current \texttt{Focus stage}, resulting in three corresponding prompts for Stage~I execution failures, Stage~II output failures, and Stage~III feasibility failures.

\begin{promptbox}{BROAD Guidance: Stage I --- Execution}{broad-execution}
\small
\raggedright

\textbf{Improvement Guidance --- Mode: \textsc{Broad}} \hfill \textit{Stage I focus: execution}

\medskip

Many cases are failing at execution---the program crashes, raises exceptions, or times out before producing output.

\medskip

Parameter-level, Strategy-level, and Algorithm-level changes are all permitted in this mode. The root cause of execution failures can sit at any level---if persistent timeouts suggest the current algorithmic paradigm is too expensive, switching to a paradigm with lower complexity is a reasonable option.

\medskip

\textbf{Please focus on:}

\begin{promptlist}

\item[1.] Eliminating syntax, import, and runtime errors so the program runs end-to-end for every case.

\item[2.] Using only the allowed libraries.

\item[3.] Keeping the program within the time limit.

\end{promptlist}

\end{promptbox}

\begin{promptbox}{BROAD Guidance: Stage II --- Valid Output}{broad-output}
\small
\raggedright

\textbf{Improvement Guidance --- Mode: \textsc{Broad}} \hfill \textit{Stage II focus: valid output}

\medskip

The program runs without crashing but does not produce valid output---the output file is empty, missing, or in the wrong format.

\medskip

Parameter-level, Strategy-level, and Algorithm-level changes are all permitted in this mode. The root cause of output failures can sit at any level.

\medskip

\textbf{Please focus on:}

\begin{promptlist}

\item[1.] Producing a non-empty, correctly formatted output file for every case, matching the format specified in the problem description.

\item[2.] Verifying the output file is created at the exact path the verifier expects.

\item[3.] Fixing edge cases where the program silently produces no output.

\end{promptlist}

\end{promptbox}

\begin{promptbox}{BROAD Guidance: Stage III --- Feasibility}{broad-feasibility}
\small
\raggedright

\textbf{Improvement Guidance --- Mode: \textsc{Broad}} \hfill \textit{Stage III focus: feasibility}

\medskip

Stage~III is the dominant failure stage. In affected cases, the solver reaches verification, but the generated solutions are rejected as infeasible. The primary objective of this round is therefore to restore solution feasibility.

\medskip

Parameter-level, Strategy-level, and Algorithm-level changes are all permitted in this mode.

\medskip

\textbf{Please focus on:}

\begin{promptlist}

\item[1.] Repairing the recurring verifier failures using the provided failure evidence.

\item[2.] Making solutions feasible first; a valid but suboptimal solution is preferable to an invalid one at this stage.

\item[3.] Fixing the underlying defect at the appropriate revision level.

\item[4.] If the current algorithmic paradigm appears mismatched with the problem structure, replacing it is a reasonable option.

\item[5.] Do not optimise for cost in this round; feasibility comes first.

\end{promptlist}

\medskip

If a few cases already pass, that is not yet evidence that the algorithm is correct---those cases may simply not exercise the faulty path. Treat this round as feasibility recovery regardless.

\end{promptbox}

\textbf{\textsc{Moderate}.} When $0.4 \leq \rho_t < 0.8$, several demonstration cases already succeed, indicating that the current solver has established partially verified behavior. DiSH uses Strategy-level revision as the default in this mode, while allowing Parameter-level changes at any time. Algorithm-level revision is reserved for failures that suggest a structural mismatch in the current algorithmic paradigm.

\begin{promptbox}{MODERATE Guidance}{moderate}
\small
\raggedright

\textbf{Improvement Guidance --- Mode: \textsc{Moderate}} \hfill \textit{\ph{SUCCESS\_RATIO} Success}

\medskip

Several demo/train cases already pass---the algorithmic skeleton has been validated on real instances. Some cases still fail; they may fail at any pipeline stage (execution, output, or verification). This round targets those remaining failures: make them pass first, then improve solution quality.

\medskip

Parameter-level, Strategy-level, and Algorithm-level changes are all permitted in this mode. This mode is NOT a conservative mode and is distinct from the 80\%+ CONSERVATIVE mode, which forbids Algorithm-level changes. Within the permitted set, prefer the smaller change scope first: make Strategy-level changes by default, and escalate to Algorithm-level changes when the failing cases' symptoms point to a structural mismatch in the current paradigm. Parameter-level changes are available at any time.

\medskip

Fixing incorrect implementation is not restricted by the granularity above: if a failing case is caused by a defect---a wrong loop bound, a wrong index, a mishandled constraint---correct it directly, even if the defect lies in code that currently passing cases also run through.

\medskip

\textbf{Please focus on:}

\begin{promptlist}

\item[1.] Use the provided failure evidence to repair the remaining failing cases.

\item[2.] After any change, the cases that passed before must still pass. Manage regression by checking this end result---not by treating any particular block of existing code as off-limits.

\item[3.] Securing feasibility on every case before attempting any cost improvement.

\end{promptlist}

\end{promptbox}

\textbf{\textsc{Conservative}.} When $0.8 \leq \rho_t < 1$, most demonstration cases already succeed, and only a small number of failures remain. DiSH limits the revision scope accordingly: Parameter-level changes are allowed, Strategy-level changes must be narrowly targeted at the remaining failures, and Algorithm-level revisions are disabled to reduce the risk of regression.

\begin{promptbox}{CONSERVATIVE Guidance}{conservative}
\small
\raggedright

\textbf{Improvement Guidance --- Mode: \textsc{Conservative}} \hfill \textit{\ph{SUCCESS\_RATIO} Success}

\medskip

Most demo/train cases are passing; only a small number remain failing. The case-by-case feedback above shows exactly which cases failed and how.

\medskip

Parameter-level changes are permitted. Strategy-level changes are permitted only when narrowly targeted at the failing cases shown above. Algorithm-level changes are not permitted in this mode---they carry a high risk of regressing the currently passing cases.

\medskip

Correctness fixes are exempt from this granularity limit (see point 5 of the system instructions).

\medskip

\textbf{Please focus on:}

\begin{promptlist}

\item[1.] Diagnosing the specific failures shown in the case feedback above, and making the smallest change that addresses them.

\item[2.] Avoiding edits to logic that has no causal link to those failing cases.

\item[3.] When uncertain whether a change might affect passing cases, prefer not to make it.

\end{promptlist}

\end{promptbox}

\textbf{\textsc{Protected}.} When $\rho_t=1$, all demonstration cases succeed, and the current solver has established fully verified behavior on the demonstration set. DiSH therefore places the solver in a protected state, where preserving the current behavior becomes the default. Strategy and Algorithm level revisions are disabled, while Parameter level changes are allowed only when they provide a clear benefit without adversely affecting any currently successful case.

\begin{promptbox}{PROTECTED Guidance}{protected}
\small
\raggedright

\textbf{Improvement Guidance --- Mode: \textsc{Protected}} \hfill \textit{100\% Success}

\medskip

All demo/train cases currently pass. The solver is in a frozen state. The default action in this mode is to return the existing solver code unchanged. This is the correct response, not a failure to follow instructions. In most rounds at 100\% success, no change should be made.

\medskip

Only Parameter-level changes are permitted in this mode. Strategy-level and Algorithm-level changes are not permitted.

\medskip

\textbf{Specifically:}

\begin{promptlist}

\item[--] Do NOT add new functions, classes, or data structures.

\item[--] Do NOT add new control flow branches or new logic paths.

\item[--] Do NOT refactor existing code, even if a cleaner version would be equivalent.

\item[--] Do NOT replace existing logic with logic you consider ``better''.

\item[--] Do NOT make speculative efficiency edits without a concrete reason.

\end{promptlist}

\medskip

A Parameter-level edit is acceptable only if you can point to a concrete, specific reason it improves the solver, AND you can confirm it cannot adversely affect any currently passing case. If either condition is unclear, return the solver unchanged.

\medskip

Correctness fixes are exempt from this restriction (see point 5 of the system prompt).

\end{promptbox}

\subsection{Revision Direction Construction}\label{app:direction}
Revision direction construction determines which failures should receive attention and what evidence should guide the next solver revision. Rather than exposing all instance-level failures to the LLM, DiSH groups failures with shared signatures and selects representative evidence from these groups, while retaining successful cases as behavioral anchors to discourage regression. Based on this common representation, we consider three direction strategies, GWRS, CRS, and GRS, which differ in how failure evidence is prioritized, selected, and organized. The resulting guidance is dynamically constructed from the evaluation records available at each refinement round.

\textbf{Group-Wise Representative Sampling (GWRS).} GWRS groups failures by shared signatures and retains one representative case from each group, together with all successful cases. The feedback first summarizes the selected representative evidence and then presents the corresponding case-level examples, reducing redundancy while preserving distinct failure patterns.

\begin{tcolorbox}[
    title={GWRS Feedback},
    colback=white,
    colframe=black!75,
    boxrule=0.6pt,
    arc=1.2mm,
    left=1.2mm,right=1.2mm,top=1mm,bottom=1mm
]
\small
\textbf{Summary.} These are the test cases and results from the previous iteration. Feedback uses group-wise representative sampling: \textbf{\{N\_SHOWN\_FAILING\}} failing case(s) shown out of \textbf{\{N\_FAILING\}} failing, with one representative per distinct failure group, plus all \textbf{\{N\_SUCCESS\}} successful case(s). Cases with the same failure-group signature as a representative are omitted this round as redundant.

\vspace{0.6em}
\textbf{Representative Failure Case.}

\begin{tcolorbox}[
    colback=gray!4,
    colframe=black!30,
    boxrule=0.4pt,
    arc=0.8mm,
    left=1mm,right=1mm,top=0.8mm,bottom=0.8mm
]
\textbf{Test Case 1:} \texttt{\{test\_name\}}\\
\textbf{Input file:} \texttt{\{content\}}\\
\textbf{Observed result:} \texttt{\{execution\_message\}}
\end{tcolorbox}

\vspace{0.5em}
\textbf{Successful Anchor.}

\begin{tcolorbox}[
    colback=gray!4,
    colframe=black!30,
    boxrule=0.4pt,
    arc=0.8mm,
    left=1mm,right=1mm,top=0.8mm,bottom=0.8mm
]
\textbf{Test Case 2:} \texttt{\{test\_name\}}\\
\textbf{Input file:} \texttt{\{content\}}\\
\textbf{Observed result:} \texttt{\{execution\_message\}}
\end{tcolorbox}
\end{tcolorbox}

\textbf{Correction Representative Sampling (CRS).} Unlike GWRS, which represents failures based only on their current signatures, CRS retains one smallest-input representative from each current failure group, together with all successful cases. Cross-round comparisons identify regressions, which are emphasized in the resulting revision guidance.

\begin{tcolorbox}[
    title={CRS Feedback},
    colback=white,
    colframe=black!75,
    boxrule=0.6pt,
    arc=1.2mm,
    left=1.2mm,right=1.2mm,top=1mm,bottom=1mm
]
\small

\textbf{Regression Alert.} The following cases were VALID in the previous iteration but FAIL now. These regressions must be repaired before optimizing new failures.

\begin{tcolorbox}[
    colback=gray!4,
    colframe=black!30,
    boxrule=0.4pt,
    arc=0.8mm,
    left=1mm,right=1mm,top=0.8mm,bottom=0.8mm
]
\textbf{\{case\_id\}:} was valid in the previous round \texttt{\{previous\_cost\}}\\
$\rightarrow$ now Stage \texttt{\{stage\}}: \texttt{\{execution\_message\}}
\end{tcolorbox}

\vspace{0.6em}
\textbf{Failure Groups.} Failing cases are grouped by a shared failure signature. The smallest-input case in each group is shown as its representative.

\vspace{0.5em}
\textbf{Representative Failure Case.}

\begin{tcolorbox}[
    colback=gray!4,
    colframe=black!30,
    boxrule=0.4pt,
    arc=0.8mm,
    left=1mm,right=1mm,top=0.8mm,bottom=0.8mm
]
\textbf{Test Case 1:} \texttt{\{test\_name\}}\\
\textbf{Input file:} \texttt{\{content\}}\\
\textbf{Observed result:} \texttt{\{execution\_message\}}
\end{tcolorbox}

\vspace{0.5em}
\textbf{Successful Anchor.}

\begin{tcolorbox}[
    colback=gray!4,
    colframe=black!30,
    boxrule=0.4pt,
    arc=0.8mm,
    left=1mm,right=1mm,top=0.8mm,bottom=0.8mm
]
\textbf{Test Case 2:} \texttt{\{test\_name\}}\\
\textbf{Status:} valid\\
\textbf{Cost:} \texttt{\{cost\}}
\end{tcolorbox}
\end{tcolorbox}

\textbf{Scope-Aware Group Representative Sampling (GRS).} GWRS retains one representative from every current failure group, while CRS further prioritizes cases that regress across refinement rounds. GRS instead performs group-level selection, retaining only the failure groups most relevant to the current round. Groups at the \texttt{Focus stage} are prioritized first, and the active revision mode determines how many high-priority groups are selected. One representative case is then retained from each selected group, together with all successful cases.

\begin{tcolorbox}[
    title={GRS Feedback },
    colback=white,
    colframe=black!75,
    boxrule=0.6pt,
    arc=1.2mm,
    left=1.2mm,right=1.2mm,top=1mm,bottom=1mm
]
\small

\textbf{Failure Groups Selected for This Round.}

Failing cases are grouped by a shared failure signature, and only the highest-priority groups are selected for repair this round. Groups at the focus stage are prioritized first. Mode \textbf{\{MODE\}}; \textbf{\{N\_SELECTED\}} group(s) selected.

\vspace{0.6em}
\textbf{Representative Failure Case.}

\begin{tcolorbox}[
    colback=gray!4,
    colframe=black!30,
    boxrule=0.4pt,
    arc=0.8mm,
    left=1mm,right=1mm,top=0.8mm,bottom=0.8mm
]
\textit{Selected failure group \{G\}; signature \texttt{\{signature\}}; covers \{N\_GROUP\} cases.}\\
\textbf{Test case:} \texttt{\{test\_name\}}\\
\textbf{Input:} \texttt{\{content\}}\\
\textbf{Observed result:} \texttt{\{execution\_message\}}
\end{tcolorbox}

\vspace{0.6em}
\textbf{Success Anchors --- already valid, must NOT regress.}

\begin{quote}
\texttt{- \{case\_id\}: cost \{cost\}}\\
\texttt{- \{case\_id\}: cost \{cost\}}
\end{quote}

\end{tcolorbox}

\subsection{Complete Refinement Example}\label{app:complete_round}
To illustrate how the preceding components are instantiated in a complete refinement round, we present a real example from Protein Sequence Design. The example combines the current harness state, representative failure feedback, successful cases, and mode-specific revision guidance. The refinement prompt is combined with the fixed system and user prompts and the current solver $h_t$. These inputs are then provided to the LLM to generate the next solver $h_{t+1}$.

\begin{promptbox}{Instantiated Refinement Prompt}{protein_round}
\small
\raggedright

\textbf{\# Feedback from Previous Iteration (Iteration 0)}

\medskip

\textbf{\#\# Demo/Train Results Summary}

\begin{promptlist}

\item[--] Total demo/train cases: 5

\item[--] Stage I failures (execution error): 0 (0.0\%)

\item[--] Stage II failures (no/invalid output): 0 (0.0\%)

\item[--] Stage III failures (verifier rejected): 3 (60.0\%)

\item[--] Successful cases: 2 (40.0\%)

\item[--] \textbf{Focus stage: Stage III}

\item[--] \textbf{Current mode: \textsc{Moderate}}

\item[--] Mode reason: SUCCESS is 40\%; the algorithmic skeleton has been validated, so strategy-level improvement is the default.

\end{promptlist}

\medskip

These are the test cases and results from the previous iteration. Feedback uses group-wise representative sampling: 1 failing case shown out of 3 failing, with one representative per distinct failure group, plus all 2 successful cases. Cases with the same failure-group signature as a representative are omitted this round as redundant.

\medskip

\textbf{\#\# Representative Failure Case}

\medskip

\textbf{Test Case: 1a4y}

\medskip

\textit{Representative of 3 demonstration cases in the same failure group this round; the other 2 are grouped and not shown.}

\medskip

\textbf{Result:} Invalid solution

\textbf{Error:} Verification failed: Generated sequence length (460) does not match expected length (1166).

\medskip

\textbf{\#\# Successful Anchors}

\medskip

\textbf{Test Case: 1aaj}

\textbf{Result:} Valid solution with cost 45.714285714285715

\medskip

\textbf{Test Case: 1aba}

\textbf{Result:} Valid solution with cost 40.229885057471265

\medskip

\textbf{\# Improvement Guidance --- Mode: \textsc{Moderate}} \hfill \textit{40\% Success}

\medskip

Several demo/train cases already pass---the algorithmic skeleton has been validated on real instances. Some cases still fail; they may fail at any pipeline stage (execution, output, or verification). This round targets those remaining failures: make them pass first, then improve solution quality.

\medskip

Parameter-level, Strategy-level, and Algorithm-level changes are all permitted in this mode. Within the permitted set, prefer the smaller change scope first: make Strategy-level changes by default, and escalate to Algorithm-level changes when the failing cases' symptoms point to a structural mismatch in the current paradigm. Parameter-level changes are available at any time.

\medskip

\textbf{Please focus on:}

\begin{promptlist}

\item[1.] Use the provided failure evidence to repair the remaining failing cases.

\item[2.] After any change, the cases that passed before must still pass.

\item[3.] Secure feasibility on every case before attempting any cost improvement.

\end{promptlist}

\end{promptbox}

\section{\textbf{Metric Definitions}}\label{app:metrics}
We employ three groups of metrics to evaluate the magnitude of solver revisions, the performance of the resulting solutions, and the dynamics of the refinement process.

\subsection{Revision Distance}
\label{app:revision_distance}

Section~\ref{subsec:overview} uses revision distance to quantify the magnitude of each solver revision. For a solver program $h$, we tokenize its Python source code and remove comments and formatting-related tokens, including indentation, line breaks, and end-of-file markers. Let $(z_1,\ldots,z_m)$ denote the resulting token sequence. We define $T(h)=\{(z_i,z_{i+1},\ldots,z_{i+4}) \mid 1\leq i\leq m-4\}$ as the set of all contiguous five-token subsequences.

\textbf{Jaccard Similarity.} For two consecutive solvers $h_t$ and $h_{t+1}$, we measure their code similarity using the Jaccard index:
\[
J(h_t,h_{t+1}) = \frac{ \left|T(h_t)\cap T(h_{t+1})\right| }{ \left|T(h_t)\cup T(h_{t+1})\right| }.
\]

\textbf{Revision Distance.} The revision distance is defined as the complement of their Jaccard similarity:
\[
\Delta(h_{t+1},h_t) = 1-J(h_t,h_{t+1}).
\]

A smaller $\Delta(h_{t+1},h_t)$ indicates greater code similarity between consecutive solvers, whereas a larger value indicates a more extensive revision. We do not normalize variable names, canonicalize abstract syntax trees, or otherwise account for semantic equivalence between programs. Thus, $\Delta(h_{t+1},h_t)$ serves as a transparent proxy for code-level revision magnitude rather than semantic distance.

\subsection{Solution Performance Metrics}\label{app:qyi}

Following the HeuriGym evaluation protocol~\citep{heurigym}, we evaluate each delivered solver on the held-out test instances using three metrics: Quality, Yield, and the Quality-Yield Index (QYI). Let $\mathcal{V}$ denote the set of test instances that pass verification, and let $N$ be the total number of test instances. For each $x\in\mathcal{V}$, the normalized solution quality relative to the reference cost is defined as
\[
q_x=\operatorname{Normalize}_p(c_x,c_x^\star).
\]

\textbf{Quality.} Quality measures the average normalized solution quality over successfully verified test instances:
\[
\mathrm{Quality} = \frac{1}{|\mathcal{V}|} \sum_{x\in\mathcal{V}} q_x.
\]

\textbf{Yield.} Yield measures the proportion of test instances that pass verification:
\[
\mathrm{Yield} = \frac{|\mathcal{V}|}{N}.
\]

\textbf{Quality-Yield Index (QYI).} QYI jointly captures solution quality and verification coverage through their harmonic mean:
\[
\mathrm{QYI} = \frac{ 2\,\mathrm{Quality}\,\mathrm{Yield} }{ \mathrm{Quality}+\mathrm{Yield} }.
\]
We set $\mathrm{QYI}=0$ when no test instance passes verification. All three metrics are computed separately for each run before aggregating the run-level results.

\subsection{Refinement Trajectory Metrics}\label{app:trajectory_metrics}

We additionally report three trajectory-level diagnostics computed on the demonstration instances. Each run contains five evaluated solver states, consisting of the shared initial solver followed by four refinement steps. Let $\rho_t$ denote the feasibility pass rate of state $t$, let $N$ denote the total number of runs, and let $\mathcal{C}$ denote the set of converged runs. For a run whose final state is fully feasible, we define the stable convergence state as $t_{\mathrm{conv}}=\min\{t:\rho_t=\rho_{t+1}=\cdots=\rho_T=1\}$. Thus, $t_{\mathrm{conv}}$ marks the beginning of the final uninterrupted fully feasible suffix rather than the first time full feasibility is reached.

\textbf{Convergence Rate.} Convergence Rate measures the fraction of runs that end in full feasibility:
\[
\mathrm{ConvergenceRate} = \frac{N_{\mathrm{conv}}}{N},
\]
where $N_{\mathrm{conv}}$ denotes the number of runs satisfying $\rho_T=1$.

\textbf{Survival Rate.} Survival Rate measures the fraction of runs that preserve full feasibility after reaching it for the first time:
\[
\mathrm{SurvivalRate} = \frac{N_{\mathrm{surv}}}{N_{\mathrm{hit}}},
\]
where $N_{\mathrm{hit}}$ denotes the number of runs that reach $\rho_t=1$ at least once, and $N_{\mathrm{surv}}$ denotes those that remain fully feasible thereafter.

\textbf{Tokens to Convergence.} Among converged runs, Tokens to Convergence measures the average cumulative number of LLM input and output tokens required to reach stable convergence:
\[
\mathrm{TokensToConv} = \frac{1}{|\mathcal{C}|} \sum_{r\in\mathcal{C}} \left( \sum_{k=1}^{t_{\mathrm{conv}}^{(r)}} c_k^{(r)} \right),
\]
where $\mathcal{C}$ denotes the set of converged runs and $|\mathcal{C}|$ its size. For each run $r\in\mathcal{C}$, $t_{\mathrm{conv}}^{(r)}$ denotes the refinement step at which stable convergence begins, and $c_k^{(r)}$ denotes the token consumption of its $k$-th refinement call.

\section{Experimental Setup}\label{app:experimental_setup}
This section details the experimental setup used to evaluate DiSH. We first describe the benchmark problems, then specify the configurations of the comparison methods. Finally, we evaluate the resulting trade-off between performance and efficiency.

\subsection{Benchmark Problems}\label{app:benchmarks}

We evaluate DiSH on nine combinatorial optimization problems spanning a diverse range of problem structures and scales. Six benchmarks are inherited from HeuriGym~\citep{heurigym}, for which we preserve the original demonstration and test splits as well as the corresponding evaluation setting to ensure direct comparability with prior work. To further assess generalization beyond the original benchmark suite, we include BPP, CVRP, and TSP from BPPLIB~\citep{bpplib}, CVRPLIB Set~X~\citep{cvrplibx}, and TSPLIB~\citep{tsplib}, respectively, with fixed and mutually disjoint demonstration and held-out test sets. This expanded benchmark suite provides a broader test of DiSH across both established and additional combinatorial optimization settings. Table~\ref{tab:benchmark_stats} summarizes the benchmark statistics, with problem sizes reported in the standard units used by each benchmark.

\begin{table}[t]
\centering
\small
\caption{Instance statistics for the nine benchmark problems.}
\label{tab:benchmark_stats}
\setlength{\tabcolsep}{4pt}

\begin{tabular}{lccccc}
\toprule
Problem & Source & \# Demo & Demo size & \# Test & Test size \\
\midrule
Operator Scheduling & HeuriGym & 5 & 3--34 operations & 19 & 28--1,500 operations \\

Intra & HeuriGym & 6 & 5--836 nodes & 22 & 771--65,335 nodes \\

Protein Design & HeuriGym & 5 & 87--1,166 residues & 19 & 46--2,958 residues \\

Pedigree & HeuriGym & 5 & 3--5 individuals & 15 & 3--350 individuals \\

Crew Pairing & HeuriGym & 5 & 37--44 flight legs & 9 & 206--6,799 flight legs \\

E-Graph & HeuriGym & 2 & 7--9 e-nodes & 21 & 1,541--121,449 e-nodes \\

BPP & BPPLIB & 5 & 100 items & 20 & 300--1,000 items \\

CVRP & CVRPLIB Set X & 5 & 119--227 customers & 20 & 297--855 customers \\

TSP & TSPLIB & 5 & 99--198 cities & 20 & 264--1,655 cities \\
\bottomrule
\end{tabular}
\end{table}

\textbf{Reference Costs and Verification.} Solution quality is evaluated against a per-instance reference cost. For the six benchmarks inherited from HeuriGym, we use the objective values of the expert-provided solutions as the per-instance reference costs~\citep{heurigym}. For BPP, CVRP, and TSP, we use the published reference objective values of the corresponding benchmark instances~\citep{bpplib,cvrplibx,tsplib}, providing a strong and consistent basis for quality comparison. For these three additional benchmarks, the evaluator also independently verifies solution feasibility and recomputes the objective value from the submitted solution rather than relying on solver-reported values. Specifically, BPP requires every item to be assigned exactly once to non-empty bins without exceeding bin capacity, CVRP requires every customer to be served exactly once while satisfying vehicle capacity, and TSP requires a tour that visits every city exactly once. The evaluator then recomputes the number of used bins, route cost under the benchmark distance convention, or complete tour length including the closing edge, respectively.

\subsection{Baseline Configurations}
\label{app:baseline_configs}
\textbf{Compared Methods.} We compare DiSH against five representative LLM-based AHD baselines: EoH~\citep{eoh}, ReEvo~\citep{reevo}, HSEvo~\citep{hsevo}, HiFo-Prompt~\citep{hifoprompt}, and RefineEvo~\citep{refineevo}. The original HeuriGym~\citep{heurigym} refinement procedure serves as an additional reference. We evaluate three DiSH variants, DiSH-GWRS, DiSH-CRS, and DiSH-GRS, corresponding to the revision direction constructions introduced in Section~\ref{subsec:direction}. All methods operate in the same complete-solver refinement setting. For the five AHD baselines, adaptations are limited to incorporating the shared initialization, complete-solver evaluation and output adaptation, minimal compatibility adjustments, and token accounting and logging. Their native search operators and method-specific search procedures remain unchanged. All AHD baselines use the same demonstration-set scalar fitness, without method-specific fitness tuning or held-out feedback.

\textbf{Shared Initialization.} For each problem, we generate three independent initial complete solvers. Within each initialization, all methods start from the same solver, denoted by $\theta_0$, so that differences in starting conditions do not confound the comparison. Each method is evaluated with three independent repeats per initialization under each backbone, yielding nine runs for every combination of problem and backbone. For the population-based AHD baselines, initialization is adapted to incorporate $\theta_0$ while preserving the subsequent search procedure. EoH, HiFo-Prompt, and RefineEvo initialize generation zero with the singleton population $\{\theta_0\}$. ReEvo and HSEvo first evaluate $\theta_0$ and then generate seed-conditioned descendants, retaining $\theta_0$ as an elite or best reference. Subsequent parent selection, reflection, crossover, mutation, and population updates remain method-specific. HeuriGym also starts from $\theta_0$ and retains its original refinement procedure. Table~\ref{tab:baseline_configs} summarizes the initialization and search configurations for the five AHD baselines.

\textbf{Search and Refinement Settings.} DiSH follows the $i=5$ setting of HeuriGym, with the shared initial solver $\theta_0$ replacing program generation in the first iteration. It then performs four refinement steps, yielding five evaluated solver states in total. For DiSH, HeuriGym, and the corresponding ablation configurations, we retain one round of refinement history in all runs. Each AHD baseline follows its own search procedure and is configured to allow at least one complete search cycle. EoH, ReEvo, and HSEvo evaluate at most eight additional solvers beyond $\theta_0$, while HiFo-Prompt and RefineEvo evaluate at most six. For ReEvo and HSEvo, \texttt{max\_fe=9} includes the evaluation of $\theta_0$. We do not enforce identical candidate or iteration counts across methods to preserve their native search behavior. Each baseline retains its original resampling and termination rules, without additional fallback or rescue candidates. The number and structure of LLM calls also vary across methods. We therefore measure search cost by the total LLM input and output tokens consumed across all calls in a completed run.

\begin{table*}[t]
\centering
\small
\caption{Initialization and search configurations for the five AHD baselines under the shared initial solver $\theta_0$. The new-solver cap excludes $\theta_0$.}
\label{tab:baseline_configs}

\setlength{\tabcolsep}{4pt}
\begin{tabularx}{\linewidth}{@{}
  >{\raggedright\arraybackslash}p{0.13\linewidth} >{\raggedright\arraybackslash}p{0.29\linewidth} >{\centering\arraybackslash}p{0.11\linewidth} >{\raggedright\arraybackslash}X@{}}
\toprule
Method & Shared-$\theta_0$ initialization & New-solver\newline cap & Key configuration \\
\midrule

EoH & Singleton population $\{\theta_0\}$, bypassing the native initialization step & $\leq 8$ & \texttt{use\_seed=True}, \texttt{pop\_size=2}, \texttt{n\_pop=4}, \texttt{max\_sample\_nums=8} \\

ReEvo & $\theta_0$ is evaluated first, followed by two seed-conditioned descendants & $\leq 8$ & \texttt{init\_pop\_size=2}, \texttt{pop\_size=2}, \texttt{mutation\_rate=0.5}, \texttt{max\_fe=9} \\

HSEvo & $\theta_0$ is evaluated first, followed by two seed-conditioned descendants & $\leq 8$ & \texttt{init\_pop\_size=2}, \texttt{pop\_size=2}, \texttt{mutation\_rate=0.5}, \texttt{hm\_size=2}, \texttt{max\_iter=1}, \texttt{max\_fe=9} \\

HiFo-Prompt & Singleton population $\{\theta_0\}$ through the native seeded initialization path & $\leq 6$ & \texttt{pop\_size=2}, \texttt{n\_pop=1}; operator E1 is applied twice, and E2, M1, M2, and M3 once each \\

RefineEvo & Singleton generation zero $\{\theta_0\}$, with no separately sampled initialization batch & $\leq 6$ & \texttt{init\_pop\_size=0}, \texttt{pop\_size=2}, \texttt{n\_pop=12}, \texttt{offspring\_per\_operator=1}, $B_{\mathrm{new}}=6$ \\

\bottomrule
\end{tabularx}
\end{table*}

\subsection{Search Cost and Efficiency}
\label{app:search_cost}

\textbf{Token Consumption.} We measure search cost by the cumulative input and output token consumption across all LLM calls in a completed run. The reported token counts are obtained directly from the LLM provider. Token consumption is first averaged over repeated runs for each combination of problem, backbone, and initialization. We then average across the three initializations and compute a macro average over the nine benchmark problems. Because LLM providers may use different tokenization schemes and accounting rules, absolute token counts are compared only within the same backbone. For relative comparisons, token consumption is normalized to DiSH-CRS under each backbone. Table~\ref{tab:search_cost} reports the resulting token consumption for each backbone, together with values normalized to DiSH-CRS.

\begin{table}[t]
\centering
\footnotesize
\caption{Cumulative LLM token consumption across backbone models.}
\label{tab:search_cost}
\begin{adjustbox}{max width=\textwidth}
\begin{tabular}{lcccc}
\toprule
Method & Qwen3.5-122B-A10B & DeepSeek-V3 & GPT-4.1-mini & Gemini-2.5-Flash \\
\midrule
EoH & 102{,}504 (3.45$\times$) & 34{,}228 (1.45$\times$) & 45{,}519 (1.94$\times$) & 37{,}037 (1.70$\times$) \\

ReEvo & 97{,}541 (3.29$\times$) & 43{,}182 (1.83$\times$) & 71{,}916 (3.07$\times$) & 48{,}055 (2.20$\times$) \\

HSEvo & 51{,}594 (1.74$\times$) & 29{,}816 (1.27$\times$) & 46{,}582 (1.99$\times$) & 39{,}197 (1.80$\times$) \\

HiFo-Prompt & 79{,}241 (2.67$\times$) & 29{,}763 (1.26$\times$) & 48{,}110 (2.05$\times$) & 30{,}048 (1.38$\times$) \\

RefineEvo & 98{,}298 (3.31$\times$) & 60{,}809 (2.58$\times$) & 68{,}114 (2.90$\times$) & 53{,}570 (2.46$\times$) \\

HeuriGym & 37{,}485 (1.26$\times$) & 28{,}075 (1.19$\times$) & 28{,}378 (1.21$\times$) & 26{,}103 (1.20$\times$) \\

\midrule
DiSH-GWRS & 30{,}133 (1.02$\times$) & 25{,}959 (1.10$\times$) & 26{,}657 (1.14$\times$) & 25{,}484 (1.17$\times$) \\

DiSH-CRS & 29{,}684 (1.00$\times$) & 23{,}536 (1.00$\times$) & 23{,}454 (1.00$\times$) & 21{,}808 (1.00$\times$) \\

DiSH-GRS & 27{,}065 (0.91$\times$) & 23{,}122 (0.98$\times$) & 23{,}478 (1.00$\times$) & 22{,}198 (1.02$\times$) \\
\bottomrule
\end{tabular}
\end{adjustbox}
\end{table}

\section{Evidence of Failure Feedback Redundancy}\label{app:failure_signatures}

\begin{figure}[t]
    \centering
    \includegraphics[width=0.7\linewidth]{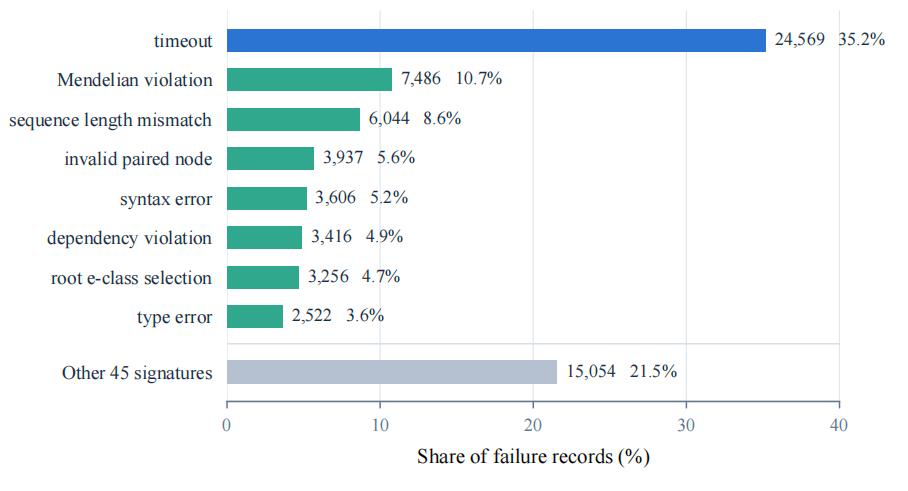}
    \caption{Distribution of normalized failure patterns across the 69,890 failure records collected from all nine benchmarks and refinement rounds. The eight most frequent patterns are shown individually, while the remaining 45 are grouped as \texttt{Other}.}
    \label{fig:failure_distribution}
\end{figure}

\begin{figure}[t]
    \centering
    \includegraphics[width=1.0\linewidth]{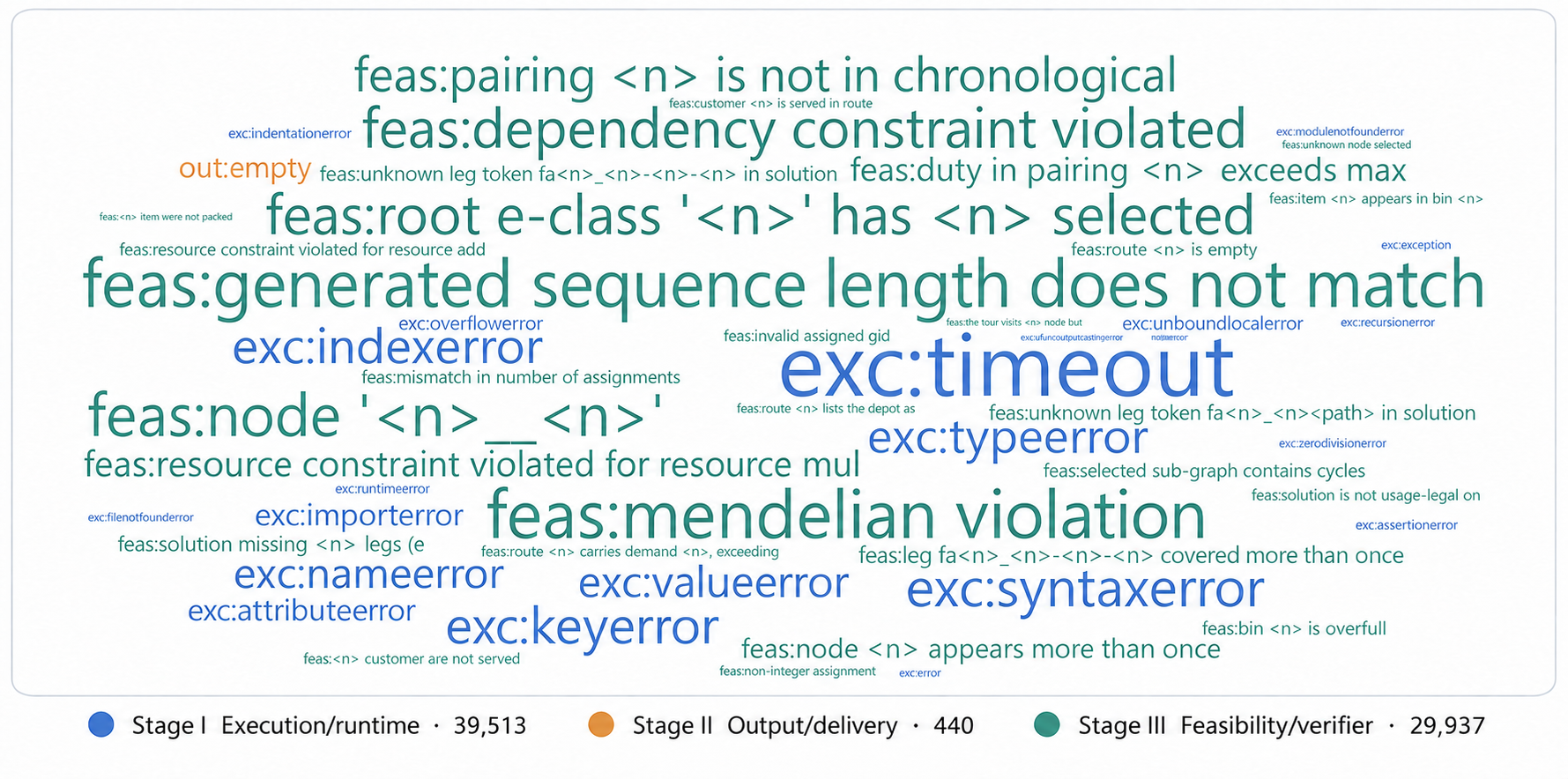}
    \caption{Visualization of the 53 normalized failure patterns. Word area scales with the square root of the observed frequency, while colors indicate the corresponding terminal evaluation stage: Stage~I (execution/runtime), Stage~II (output/delivery), and Stage~III (feasibility/verifier).}
    \label{fig:failure_cloud}
\end{figure}

To examine redundancy in the failure feedback encountered during solver refinement, we aggregate recurring failure patterns across all nine benchmarks and refinement rounds. The 69,890 recorded failures collapse into only 53 distinct normalized failure patterns, revealing substantial redundancy in instance-level feedback. This observation motivates the use of representative failure evidence during solver refinement. Figure~\ref{fig:failure_distribution} shows the frequency distribution of these patterns, while Figure~\ref{fig:failure_cloud} provides a complementary view of their recurrence and associated evaluation stages.

\section{Additional Performance Analysis}\label{app:additional_results}
This section provides additional analyses of DiSH beyond the aggregate QYI results. We first examine the contributions of Quality and Yield, including a closer analysis of Pedigree. We then study performance across problems and backbone models, as well as across shared initializations.

\subsection{Quality--Yield Analysis}\label{app:qy_decomposition}

To better understand the source of the QYI improvements, we examine Quality and Yield separately at both the problem and backbone levels. All three DiSH variants improve both Quality and Yield, indicating that their QYI gains do not come at the expense of either component. Quality measures average solution quality over successfully verified test instances, while Yield measures the proportion passing verification. Tables~\ref{tab:qy_problem} and~\ref{tab:qy_backbone} report problem- and backbone-level results as Quality / Yield (\%). DiSH-CRS achieves average Quality and Yield of 80.7\% and 84.5\%, corresponding to relative improvements of 4.9\% and 23.2\% over the respective strongest non-DiSH averages, EoH's 76.9\% and HeuriGym's 68.6\%. Even the lowest averages across the three variants, 79.2\% for Quality and 81.6\% for Yield, exceed these non-DiSH maxima. Across backbones, all three variants consistently outperform every non-DiSH method in Yield, although their Quality advantages are not universal. Problem-level patterns also vary: on Pedigree, all three variants have lower Quality than EoH but not consistently lower Yield, motivating further analysis of this benchmark.

\begin{table*}[t]
\centering
\small
\caption{Problem-level Quality and Yield decomposition. Each entry reports Quality/Yield (\%).}
\label{tab:qy_problem}
\begin{adjustbox}{max width=\textwidth}
\begin{tabular}{lcccccccccc}
\toprule
Method & Operator & Intra & Protein & Pedigree & Crew & E-Graph & BPP & CVRP & TSP & Avg. \\
\midrule

EoH & 90.9 / 84.8 & 71.6 / 72.9 & 73.7 / 61.8 & 81.2 / 59.4 & 64.6 / 58.6 & 41.9 / 43.7 & 99.4 / 100.0 & 88.1 / 86.7 & 80.7 / 40.4 & 76.9 / 67.6 \\

ReEvo & 71.3 / 66.1 & 63.0 / 54.7 & 67.1 / 53.4 & 61.3 / 62.4 & 38.0 / 19.8 & 26.9 / 27.9 & 99.4 / 93.3 & 83.2 / 75.4 & 57.4 / 32.9 & 63.1 / 54.0 \\

HSEvo & 72.9 / 64.3 & 57.4 / 59.3 & 68.2 / 65.8 & 56.5 / 43.0 & 44.5 / 23.5 & 40.2 / 41.7 & 99.4 / 99.6 & 81.0 / 78.3 & 57.8 / 29.3 & 64.2 / 56.1 \\

HiFo-Prompt & 85.9 / 78.5 & 62.5 / 62.0 & 73.7 / 59.6 & 68.1 / 74.1 & 46.1 / 29.0 & 35.1 / 28.3 & 99.6 / 97.9 & 78.3 / 74.7 & 52.4 / 32.1 & 66.9 / 59.6 \\

RefineEvo & 74.4 / 64.8 & 74.1 / 55.1 & 63.6 / 48.8 & 73.7 / 57.6 & 52.6 / 35.2 & 40.4 / 45.2 & 91.1 / 85.0 & 77.4 / 76.0 & 63.0 / 29.2 & 67.8 / 55.2 \\

HeuriGym & 91.2 / 93.1 & 70.1 / 71.0 & 69.4 / 82.2 & 69.0 / 41.1 & 65.7 / 40.1 & 50.0 / 55.0 & 99.4 / 100.0 & 91.5 / 83.3 & 65.5 / 51.8 & 74.6 / 68.6 \\

\midrule

DiSH-GWRS & 88.4 / 100.0 & 74.2 / 69.9 & 79.8 / 84.6 & 71.5 / 59.4 & 77.1 / 75.6 & 54.6 / 58.9 & 99.4 / 100.0 & 89.3 / 97.6 & 85.9 / 88.2 & 80.0 / 81.6 \\

DiSH-CRS & 88.1 / 100.0 & 71.0 / 73.5 & 80.4 / 84.1 & 74.4 / 58.5 & 66.6 / 79.0 & 69.3 / 78.2 & 99.4 / 100.0 & 90.2 / 97.6 & 87.2 / 89.2 & 80.7 / 84.5 \\

DiSH-GRS & 86.3 / 97.2 & 69.4 / 77.9 & 81.1 / 81.6 & 67.2 / 64.8 & 71.9 / 81.2 & 61.6 / 66.1 & 99.4 / 100.0 & 88.8 / 97.8 & 86.9 / 84.4 & 79.2 / 83.5 \\

\bottomrule
\end{tabular}
\end{adjustbox}
\end{table*}

\begin{table*}[t]
\centering
\small
\caption{Backbone-level Quality and Yield decomposition. Each entry reports Quality/Yield (\%).}
\label{tab:qy_backbone}
\begin{adjustbox}{max width=\textwidth}
\begin{tabular}{lccccc}
\toprule
Method & Qwen3.5-122B-A10B & DeepSeek-V3 & GPT-4.1-mini & Gemini-2.5-Flash & Avg. \\
\midrule

EoH & 77.2 / 71.4 & 69.1 / 56.8 & 75.6 / 69.7 & 85.8 / 72.4 & 76.9 / 67.6 \\

ReEvo & 56.3 / 45.9 & 53.5 / 43.7 & 64.2 / 56.9 & 78.4 / 69.3 & 63.1 / 54.0 \\

HSEvo & 58.2 / 52.4 & 54.7 / 44.6 & 64.3 / 55.3 & 79.6 / 72.0 & 64.2 / 56.1 \\

HiFo-Prompt & 69.4 / 64.0 & 55.5 / 45.0 & 60.1 / 60.9 & 82.5 / 68.4 & 66.9 / 59.6 \\

RefineEvo & 53.3 / 40.8 & 57.7 / 38.2 & 79.1 / 66.7 & 81.2 / 75.1 & 67.8 / 55.2 \\

HeuriGym & 71.2 / 68.1 & 73.5 / 62.5 & 70.4 / 67.9 & 83.4 / 76.0 & 74.6 / 68.6 \\

\midrule

DiSH-GWRS & 82.5 / 80.0 & 78.2 / 80.0 & 76.5 / 82.2 & 82.8 / 84.3 & 80.0 / 81.6 \\

DiSH-CRS & 79.9 / 80.8 & 81.2 / 83.9 & 80.1 / 88.6 & 81.7 / 84.5 & 80.7 / 84.5 \\

DiSH-GRS & 80.4 / 84.1 & 76.0 / 81.5 & 80.5 / 85.8 & 79.8 / 82.3 & 79.2 / 83.5 \\

\bottomrule
\end{tabular}
\end{adjustbox}
\end{table*}

\textbf{Analysis of the Pedigree Results.}\label{app:pedigree_qy} DiSH achieves the highest average QYI on all evaluated benchmarks except Pedigree. To better understand this exception, we compare the three DiSH variants with the baseline methods considered in this study. We use EoH as the main reference because it achieves the highest average QYI among the five AHD baselines (Table~\ref{tab:problem_results}). Table~\ref{tab:pedigree_qy} reports the corresponding Quality and Yield results. The results suggest that, relative to EoH, the lower QYI on Pedigree is mainly attributable to lower Quality rather than consistently lower Yield. EoH achieves an average Quality of 81.2\% and a Yield of 59.4\%. DiSH-GWRS reaches the same Yield but with lower Quality, while DiSH-CRS achieves a similar Yield with a smaller reduction in Quality. Among the three DiSH variants, DiSH-GRS obtains the highest Yield at 64.8\%, but its average Quality decreases to 67.2\%. Overall, the three variants exhibit different balances between Quality and Yield on Pedigree, indicating that higher verification coverage does not necessarily lead to better solution Quality.

\begin{table*}[t]
\centering
\small
\caption{Quality and Yield on Pedigree for EoH and the three DiSH variants. Each entry reports Quality/Yield (\%).}
\label{tab:pedigree_qy}
\begin{adjustbox}{max width=\textwidth}
\begin{tabular}{lccccc}
\toprule
Method & Qwen3.5-122B-A10B & DeepSeek-V3 & GPT-4.1-mini & Gemini-2.5-Flash & Avg. \\
\midrule

EoH & 84.0 / 67.4 & 83.3 / 47.4 & 66.1 / 67.4 & 91.5 / 55.6 & 81.2 / 59.4 \\

\midrule

DiSH-GWRS & 75.1 / 48.9 & 68.4 / 67.4 & 62.5 / 61.5 & 79.9 / 60.0 & 71.5 / 59.4 \\

DiSH-CRS & 75.0 / 60.7 & 63.8 / 57.0 & 67.3 / 64.4 & 91.5 / 51.9 & 74.4 / 58.5 \\

DiSH-GRS & 62.1 / 65.2 & 60.6 / 63.7 & 76.1 / 65.2 & 69.9 / 65.2 & 67.2 / 64.8 \\

\bottomrule
\end{tabular}
\end{adjustbox}
\end{table*}

\subsection{Performance across Problems and Backbones}\label{app:problem_backbone}

\begin{figure*}[t]
    \centering

    \begin{minipage}[h]{0.32\textwidth}
        \centering
        \includegraphics[width=\linewidth]{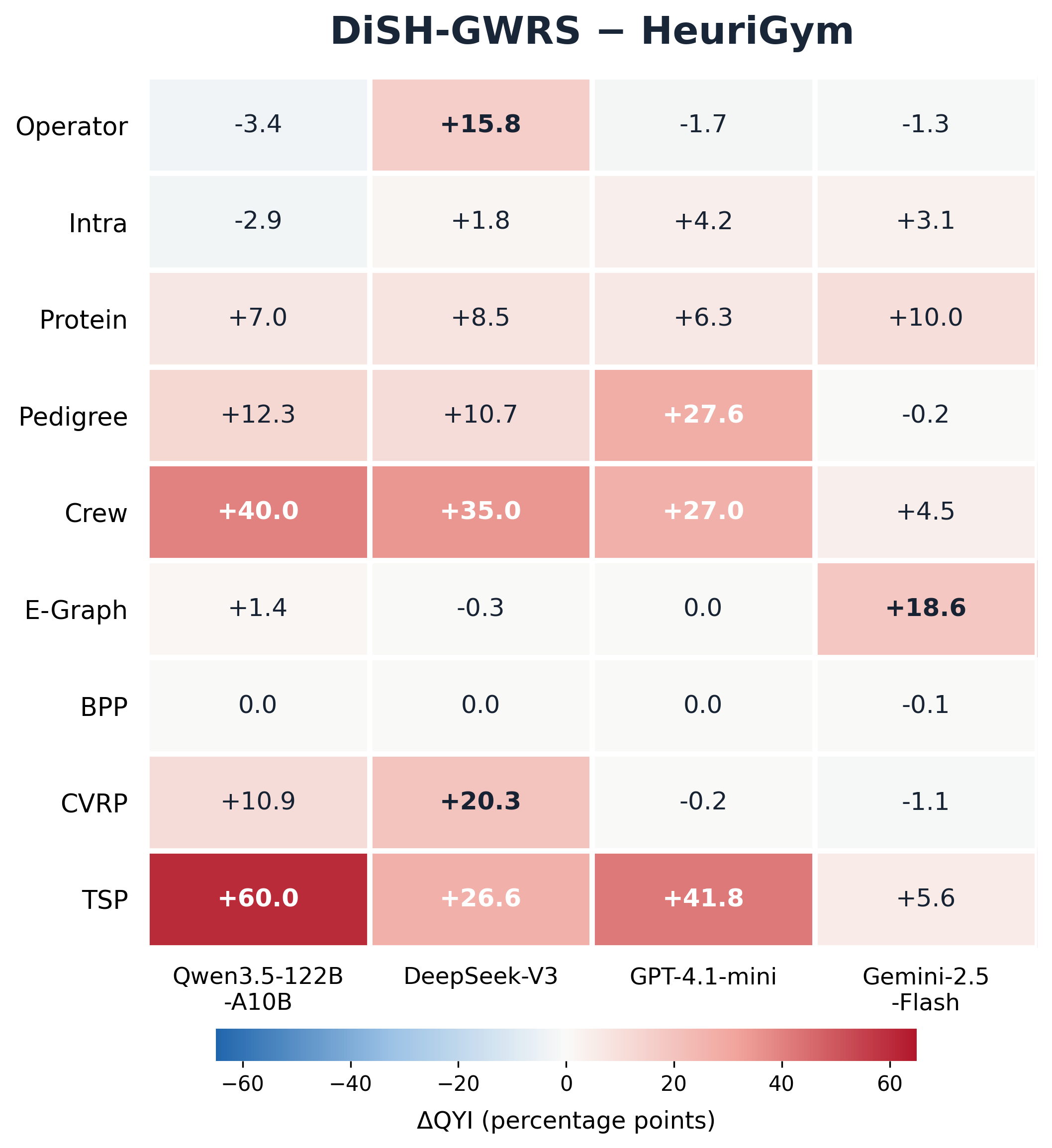}
        \vspace{-1.5mm}
        {\small (a) DiSH-GWRS}
    \end{minipage}
    \hfill
    \begin{minipage}[h]{0.32\textwidth}
        \centering
        \includegraphics[width=\linewidth]{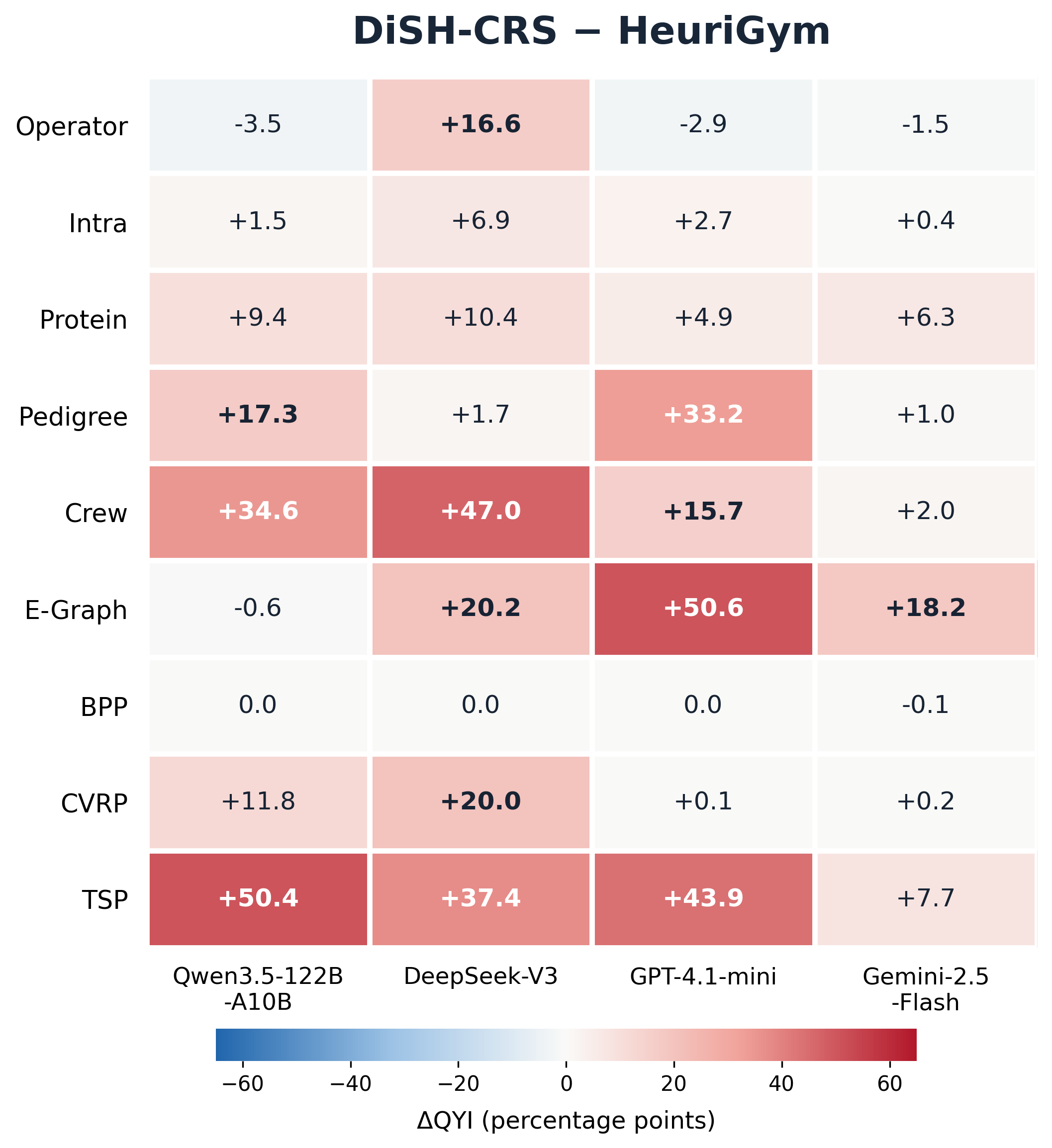}
        \vspace{-1.5mm}
        {\small (b) DiSH-CRS}
    \end{minipage}
    \hfill
    \begin{minipage}[h]{0.32\textwidth}
        \centering
        \includegraphics[width=\linewidth]{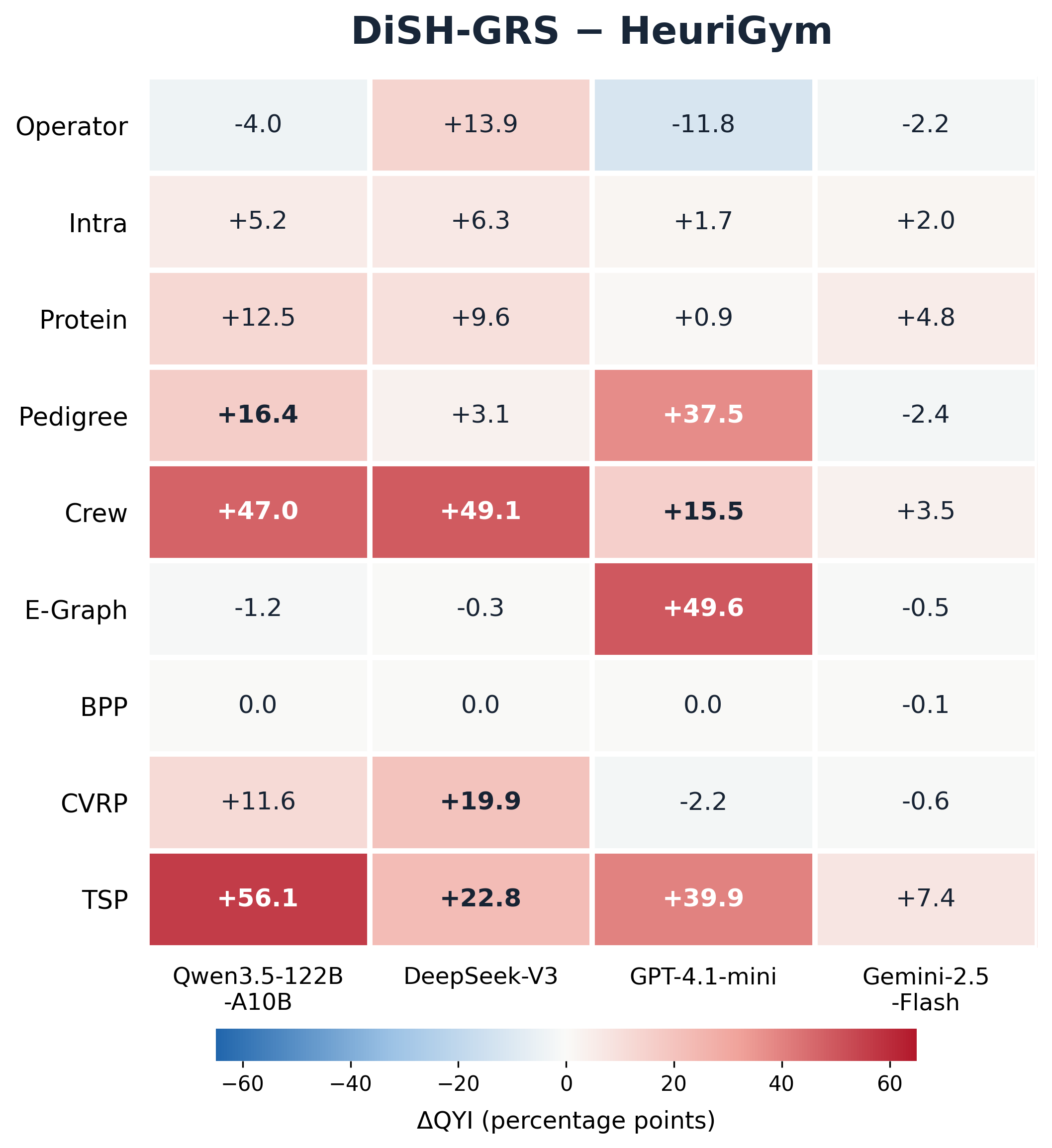}
        \vspace{-1.5mm}
        {\small (c) DiSH-GRS}
    \end{minipage}
    \caption{QYI differences between each DiSH variant and HeuriGym across problems and backbone models.}
    \label{fig:problem_backbone_heatmaps}
\end{figure*}

To examine whether the aggregate gains are concentrated in particular problems or backbone models, we further compare each DiSH variant with HeuriGym across all combinations of problems and backbones. We quantify the improvement for each problem $p$ and backbone $b$ by its QYI gain over HeuriGym: \[ \Delta \mathrm{QYI}_{p,b} = \mathrm{QYI}_{\mathrm{DiSH},p,b} - \mathrm{QYI}_{\mathrm{HeuriGym},p,b}. \] Figure~\ref{fig:problem_backbone_heatmaps} visualizes the resulting 36 differences for each DiSH variant. Each cell reports the QYI difference in percentage points, with a shared color scale centered at zero.

The gains are distributed across multiple problems and backbones rather than being driven by a single setting. The largest improvements occur on TSP, Crew Pairing, and E-Graph, with gains exceeding 40 QYI points for several backbone models. Protein Design also exhibits consistently positive gains across all four backbones and all three DiSH variants. In contrast, BPP remains close to zero across nearly all settings, consistent with its near-saturated performance in the main results. Overall, these results indicate that the aggregate gains of DiSH arise from improvements across a broad range of problem--backbone settings rather than from a small number of favorable cases.

\subsection{Performance across Initializations}
\label{app:paired}
Because all methods start from the same initial complete solver within each initialization, their performance can be compared directly at the initialization level. For each backbone, we compare DiSH against the strongest non-DiSH method according to its overall average QYI under that backbone. Since the comparison method is selected from the observed aggregate results and the twelve comparisons consist of four backbones with three shared initializations each, we report paired QYI differences and win counts descriptively rather than performing significance tests. Summary statistics are computed from the unrounded values. Table~\ref{tab:paired_detail} reports the individual paired comparisons for DiSH-CRS, while Table~\ref{tab:paired_summary} summarizes the corresponding results for all three DiSH variants.

\begin{table*}[t]
\centering
\small
\caption{Initialization-level paired comparisons between DiSH-CRS and the strongest non-DiSH method under each backbone. $\Delta$QYI denotes DiSH-CRS minus the comparison method under the same initialization.}
\label{tab:paired_detail}
\begin{adjustbox}{max width=\textwidth}
\begin{tabular}{llcccc}
\toprule
Backbone & Baseline & Init. & DiSH-CRS QYI (\%) & Baseline QYI (\%) & $\Delta$QYI (\%) \\
\midrule

Qwen3.5-122B-A10B & EoH & 1 & 84.5 & 69.1 & $+15.4$ \\

Qwen3.5-122B-A10B & EoH & 2 & 78.7 & 66.1 & $+12.6$ \\

Qwen3.5-122B-A10B & EoH & 3 & 69.6 & 67.5 & $+2.1$ \\

\midrule

DeepSeek-V3 & HeuriGym & 1 & 87.3 & 68.9 & $+18.4$ \\

DeepSeek-V3 & HeuriGym & 2 & 85.8 & 71.1 & $+14.7$ \\

DeepSeek-V3 & HeuriGym & 3 & 68.6 & 48.3 & $+20.3$ \\

\midrule

GPT-4.1-mini & EoH & 1 & 88.0 & 64.5 & $+23.5$ \\

GPT-4.1-mini & EoH & 2 & 78.3 & 68.1 & $+10.2$ \\

GPT-4.1-mini & EoH & 3 & 81.8 & 73.8 & $+8.0$ \\

\midrule

Gemini-2.5-Flash & HeuriGym & 1 & 87.7 & 82.8 & $+4.9$ \\

Gemini-2.5-Flash & HeuriGym & 2 & 72.4 & 65.2 & $+7.2$ \\

Gemini-2.5-Flash & HeuriGym & 3 & 81.1 & 81.8 & $-0.7$ \\

\bottomrule
\end{tabular}
\end{adjustbox}
\end{table*}

\begin{table*}[t]
\centering
\small
\caption{Summary of initialization-level paired QYI differences for the three DiSH variants. Each entry reports the mean $\pm$ standard deviation in percentage points over the three shared initializations, followed by the number of initializations on which the DiSH variant attains the higher QYI.}
\label{tab:paired_summary}
\setlength{\tabcolsep}{4pt}
\begin{adjustbox}{max width=\textwidth}
\begin{tabular}{lcccc}
\toprule
Variant & vs.\ EoH (Qwen) & vs.\ HeuriGym (DeepSeek) & vs.\ EoH (GPT) & vs.\ HeuriGym (Gemini) \\
\midrule

DiSH-GWRS & $+10.5 \pm 8.5$, 3/3 & $+13.1 \pm 4.3$, 3/3 & $+9.1 \pm 6.3$, 3/3 & $+4.3 \pm 7.1$, 2/3 \\

DiSH-CRS & $+10.0 \pm 7.0$, 3/3 & $+17.8 \pm 2.9$, 3/3 & $+13.9 \pm 8.4$, 3/3 & $+3.8 \pm 4.1$, 2/3 \\

DiSH-GRS & $+12.5 \pm 6.6$, 3/3 & $+13.8 \pm 2.9$, 3/3 & $+12.0 \pm 12.7$, 3/3 & $+1.3 \pm 8.3$, 2/3 \\

\bottomrule
\end{tabular}
\end{adjustbox}
\end{table*}

The paired differences are positive on average for all three DiSH variants under every backbone, and each variant achieves a higher QYI in 11 of the 12 initialization-level comparisons. For DiSH-CRS, all three paired differences are positive under Qwen3.5-122B-A10B, DeepSeek-V3, and GPT-4.1-mini. Under Gemini-2.5-Flash, two of the three differences are positive, while the remaining difference is small ($-0.7$ percentage points). Overall, these results show that the aggregate gains persist across shared initializations rather than being driven by a single favorable initialization.

\section{Limitations and Future Work}

\paragraph{Limitations.} The main limitation of \methodname{} is that its control mechanism is designed primarily around feasibility, with less explicit control over solution quality. As the feasibility pass rate increases, \methodname{} progressively narrows the revision scope to preserve verified behavior and reduce regressions. Once full feasibility is reached, revisions therefore become increasingly conservative. This improves refinement stability, but may also limit opportunities for further quality improvement. In addition, the current scope thresholds and direction strategies are manually specified. These rules are simple and interpretable, but may not adapt equally well to different problems, backbone LLMs, or refinement stages. Overall, \methodname{} provides effective control for stable feasibility, while quality-aware and more adaptive control remain open directions.

\paragraph{Future Work.} In the future, we can extend \methodname{} toward more flexible control of both feasibility and solution quality. One direction is to continue improving solution quality after feasibility has been established, while explicitly protecting verified feasible behavior. Another direction is to learn the direction and scope policies from refinement trajectories. This could replace fixed thresholds and manually designed rules with controls that adapt to the current solver state, problem structure, and backbone LLM. A further extension is to combine such control with population-based automatic heuristic design, where multiple solver candidates are refined in parallel and can provide complementary search information. Together, these directions could make revision control more adaptive while supporting continued improvement beyond feasibility.

\end{document}